\documentclass[hidelinks,onefignum,onetabnum]{siamart250211}

\ifpdf
\hypersetup{
  pdftitle={Controlling the energy jump of multistable structures using shape optimization},
  pdfauthor={Arselane Hadj Slimane, Patrick E.~Farrell, Alberto Paganini, \`Alex Ferrer},
  colorlinks,
  linkcolor={red!50!black},
  citecolor={blue!50!black},
  urlcolor={blue!80!black}
}
\fi

\usepackage{microtype}      
\usepackage{setspace}
\usepackage{hyperref}
\usepackage[T1]{fontenc}    % Font encoding
\usepackage[utf8]{inputenc} % Input encoding (UTF-8)
\usepackage{lmodern}        % Latin Modern font
\usepackage{rotating}
\usepackage{algorithm}
\usepackage{algorithmic}
\usepackage{todonotes}
\usepackage{cprotect}

\usepackage{tikz}
\usepackage{subcaption}

\usepackage[a4paper, margin=1in]{geometry} % Adjust margins and page size

\usepackage[english]{babel}  % Language support
\usepackage{csquotes}        % Context-sensitive quotation marks
\usepackage{cite}

\usepackage{xcolor}    % Colors
\usepackage{graphicx}  % Include graphics
\usepackage{float}     % Improved figure placement
\usepackage{caption}   % Customize captions for figures and tables
\usepackage{subcaption} % Subfigures and subcaptions
\usepackage{amssymb}    % More math symbols
\usepackage{amsfonts}   % Additional math fonts
\usepackage{mathtools}  % Extension to amsmath
\usepackage{mathtools}

\usepackage{booktabs}   % Professional-quality tables
\usepackage{array}      % Advanced table options
\usepackage{multirow}   % Merge cells in tables
\usepackage{colortbl}   % Add colors to table cells
\usepackage{longtable}  % Tables that span multiple pages

\usepackage{enumitem}    % Customization of lists (itemize, enumerate)
\usepackage{hyperref}    % Add hyperlinks to references and TOC
\usepackage{cleveref}    % Automatic reference formatting

\usepackage{listings}    % Code listings

\usepackage{setspace}    % Adjust line spacing
\usepackage{tikz}        % Create graphics programmatically
\usepackage{siunitx}     % Formatting for units and numbers
\usepackage{etoolbox}    % Programming tools for advanced macros
\usepackage{verbatim}    % Verbatim input and comments

\newcommand{\R}{\mathbb{R}}
\newcommand{\E}{\mathcal{E}}
\newcommand{\F}{{F}}
\newcommand{\U}{\mathcal{U}}
\newcommand{\N}{\mathbb{N}}
\newcommand{\dx}{\,\mathrm{d}x}

\newcommand{\I}{{I}}
\newcommand{\bH}{{H}}
\newcommand{\tr}{\operatorname{tr}}
\newcommand{\dd}{\,\mathrm{d}}

\title{Controlling the energy jump of multistable structures using shape optimization
\thanks{\funding{
AHS was supported by
the Oriel College Research Fund.
PEF was supported by 
the Engineering and Physical Sciences Research Council [grant number EP/W026163/1],
the Science and Technology Facilities Council [grant number UKRI/ST/B000495/1],
the Donatio Universitatis Carolinae Chair ``Mathematical modelling of multicomponent systems'',
the UKRI Digital Research Infrastructure Programme through the Science and Technology Facilities Council's Computational Science Centre for Research Communities (CoSeC),
and the Swedish Research Council under grant number~Z2021-06594 while in residence at Institut Mittag-Leffler in Djursholm, Sweden.
AF was supported by the Catalan Government through a Serra Húnter Fellowship, and by the TOMAT (PID2023-153213NAI00) and FLEXION (PCI2024-155060-2) Spanish projects.
For the purpose of open access, the authors have applied a CC BY public copyright 
licence to any author accepted manuscript arising from this submission.
No new data were generated or analysed in this work.
}}}

\author{Arselane Hadj Slimane \thanks{CMAP, CNRS, École polytechnique, Institut Polytechnique de Paris, Palaiseau, France  (\email{arselane.hadj-slimane@polytechnique.edu}).}
\and
Patrick E.~Farrell\thanks{Mathematical Institute, University of Oxford, UK 
and Mathematical Institute, Faculty of Mathematics and Physics, Charles 
University, Czechia
 (\email{patrick.farrell@maths.ox.ac.uk}).}
\and
Alberto Paganini\thanks{School of Computing and Mathematical Sciences,
University of Leicester, UK (\email{a.paganini@leicester.ac.uk}).}
\and
\`Alex Ferrer\thanks{Physics Department, Universitat Politècnica de Catalunya, Barcelona, Spain and International Centre for Numerical Methods in Engineering (CIMNE), Barcelona, Spain. (\email{alex.ferrer@upc.edu}).}}

\begin{document}
\maketitle

\begin{abstract}
Multistable systems admit distinct stable equilibria for the same loading and boundary conditions. Shape optimization of multistable structures offers a promising route to engineer the mechanical response of advanced materials and metamaterials. In this work, we develop a formulation and algorithm for controlling the energy jump associated with snap-through behaviour in hyperelastic metamaterials by optimizing the shape of the domain. The key challenge is that classical shape optimization theory assumes a single-valued domain-to-solution map, which does not hold in the multistable setting. We address this by extending the theory to settings with non-unique solutions, treating the scaling of the domain deformation as a continuation parameter and applying the implicit function theorem to give sufficient conditions for the continued existence of multiple solution branches as the shape varies. This theoretical foundation underpins a practical algorithm that targets a prescribed energy jump between stable states, enabling the systematic design of structures whose snap-through response can be tuned on demand.
%Algorithms and analysis for shape optimization problems constrained by partial differential equations (PDEs) typically start by assuming that the domain-to-solution map is single valued, i.e.~that the system of PDEs has a unique solution on any domain. What happens if this assumption does not hold? In this work we investigate shape optimization without this assumption, seeking to control the relations between multiple solutions of the same PDE for the same boundary conditions on the same domain. We consider the energy dissipation associated with snap-through behaviour in hyperelastic metamaterials, and design a formulation and algorithm to control the energy jump required to snap from one stable solution to another. We present an initial analysis, considering the domain as a bifurcation parameter and developing a generalised implicit function theorem to prove the continued existence of solution branches as we vary the shape.
\end{abstract}

\begin{keywords}
  Shape optimization, multistability, hyperelasticity
\end{keywords}

\begin{AMS}
  49Q10, 74B20, 65N30, 74S05
\end{AMS}

\section{Introduction}

Shape optimization lies at the intersection of partial differential equations (PDEs), optimization, physics, and engineering.
In most applications, shape optimization aims to choose a domain to optimize a cost functional subject to physical constraints, often expressed as PDEs. This field can be viewed as a branch of optimal control or PDE-constrained optimization, where the control variable is the shape of the domain on which the PDE is defined \cite{delfour2011shapes, allaire, NewTrendsShapeOpt}.
%Many engineering problems can be formulated as shape optimization tasks. In aerospace engineering, for instance, it is used to optimize airplane wing profiles for minimal drag and gas turbine disks for structural integrity. In the automotive sector, applications include the design of rubber bushings using differential evolution algorithms to meet target stiffness curves, and the aerodynamic shaping of car bodies to reduce pressure drag. Civil and mechanical engineering also leverage shape optimization for the design of arch dams, ship hull hydrodynamics, and 3D-printed composites with optimized fiber paths for maximum stiffness. Moreover, components such as control valves and dovetail joints can be improved by refining their geometries through gradient-based methods and adjoint simulations.

In this work, we focus on shape optimization problems constrained by PDEs that admit multiple solutions.
Specifically, we apply shape optimization to control the energy jump between pairs of distinct solutions 
to the same PDE constraint.
We build on prior work \cite{defcon2}, where we used shape optimization to control the bifurcation point at which multiple solutions emerge.
%In \cite{defcon2} the authors control a bifurcation point of non-linear PDEs with shape optimization by controlling the associated solution to the Moore--Spence system that identifies bifurcation points. 
%In \cite{delfour_minmax}, the authors established the shape regularity of the absolute minimal value of certain non-convex energy functionals irrespective of the number of absolute minima. 
As a motivating application, we consider the shape optimization of neo-Hookean hyperelastic metamaterial structures. However, our approach is easily adaptable to other applications.

\begin{figure}[htb!]
\centering
\includegraphics[scale=0.2]{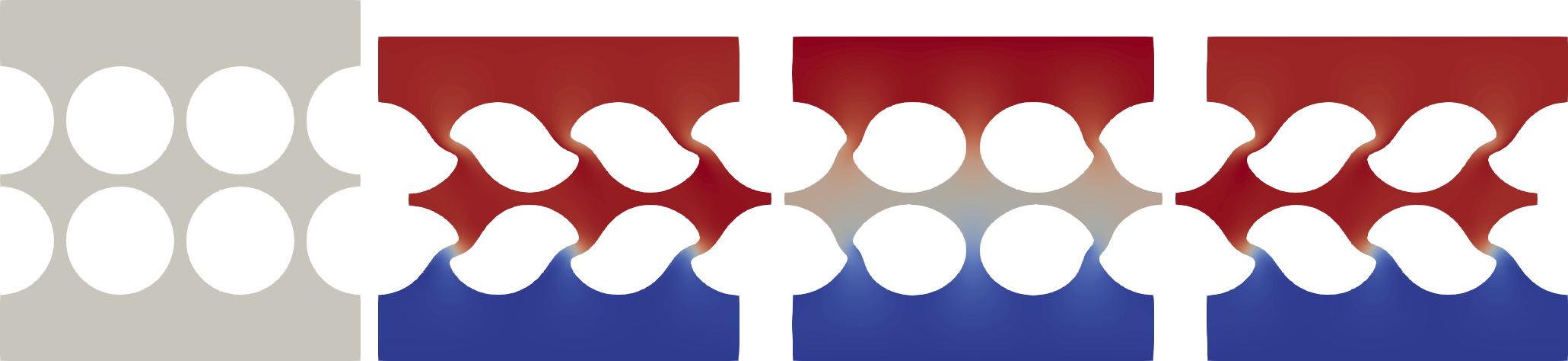}
\caption{Initial metamaterial domain $\Omega_0$ with two holes (left), and deformed domain under three distinct equilibria for the same boundary conditions:
$u_1(\Omega_0)$ (stable), $u_2(\Omega_0)$ (unstable), and $u_3(\Omega_0)$ (stable). Dirichlet conditions are imposed on the top and bottom, and natural boundary conditions on all other boundaries.}
\label{fig:Omega00}
\end{figure}

Metamaterials are engineered materials with unusual properties not found in natural materials. In particular, 
we consider elastic metamaterials that exhibit different stable solutions, %   son ratio in compression
 as shown in \Cref{fig:Omega00}.   
These materials have attracted substantial interest for their potential in engineering design for protecting against shocks \cite{metamat1,metamat2,metamat3,metamat4}. Metamaterials also exhibit a very rich bifurcation structure, with many physically achievable stable solutions \cite{defcon-metamat}. This bifurcation structure can be exploited in design. For example,
a metamaterial could dissipate impact energy by swapping branches, that is, by transitioning from one stable configuration
to another \cite{metamat-snap1,metamat-snap2,metamat-snap3}.
Controlling this energy dissipation is desirable to achieve optimal protection from impact shocks \cite{metamat-snap1,metamat-snap2,metamat-snap3}. Shape optimization has previously been applied to tune stress-strain responses and the (negative) Poisson's ratio of these materials \cite{ShapeOpt-metamat}; in this work we instead aim to control the energy required to snap from one branch to another, rendering the optimization problem qualitatively different.

Most shape optimization methods do not need to address the difficulties introduced by PDE constraints that admit multiple solutions.
Instead, it is commonly assumed that the solution to the PDE constraint is unique (possibly by introducing additional constraints). In turn, this induces a unique, single-valued shape-to-solution map, which enables the use of derivative-based optimization algorithms.
However, many physical problems exhibit a rich landscape of multiple solutions. For example, PDE constraints stemming from
nonconvex energy minimization problems can support
multiple solutions in the form of local minima, local maxima, and saddle points. In these instances, correctly accounting for the multiple solutions 
afforded by the PDE is essential.

The core tool in the bifurcation analysis of problems with multiple solutions is the implicit function theorem
(IFT) \cite[\S 7.13]{ciarlet2013}, see also \cite{IFT1,IFT2,IFT3}.
The IFT gives sufficient but not necessary conditions ensuring that solution branches can be continued locally as a function of a parameter arising in the equations. In classical bifurcation analysis, this parameter is generally a real scalar. By parameterizing the shape deformations in terms of such a parameter, we devise conditions under which the IFT guarantees the local existence of branches of solutions as the domain shape is varied. Building on this insight, we develop a shape optimization algorithm that carefully tracks these solution branches throughout the shape optimization process.

The remainder of this article is organized as follows. In \Cref{sec:problem_statement}, we formulate the optimization problem under consideration and highlight that its PDE constraint supports multiple solutions. In \Cref{sec:shapetheory}, we
propose a tailored shape optimization framework. In \Cref{sec:numerical_results}, we detail several numerical experiments.
Finally, in \Cref{sec:conclusion}, we draw the main conclusions and discuss open questions.

%As a class of initial domain shapes we focus on finite size elastomeric sheets with an embedded array of pores subjected to uniaxial tension, compression, and shear. We consider an objective function of the form:
%\begin{equation}
%\min_{\Omega \in \mathcal{U}_{ad}} J(\Omega)= \sum_{\substack{i\in I, j\in J \\ F(\Omega, u_k) = 0}} |\Delta\E(\Omega,u_i,u_j)-r\Delta\E(\Omega_0,u_i^{\Omega_0},u_j^{\Omega_0})|
%\end{equation}
%where \(F(\Omega,u)=0\) represents the PDE-constraints, in our case we consider a neo-Hookean hyperelasticity law.  Then, \(\{u_i\}_{i\in I}\) and \(\{u_j\}_{j\in J}\) are two particular subsets of solutions to the PDE, and \(\Delta\E(\Omega,u_i,u_j)\) is the energy gap between two solutions \(u_i\) and \(u_j\). \(\Omega_0\) is then the reference or initial domain and \(\{u_i^{\Omega_0})\}_{i\in I}\) and \(\{u_j^{\Omega_0})\}_{j\in J}\) are the two sets of associated initial solutions. The choice of these the sets of indeces \(I\) and \(J\) is dependent on the practical application considered, and \(r\in(0,1)\) is a chosen ratio. We aim therefore at controlling the energy jumps between different branches of solutions to a percentage \( r\) of the initial energy gap.
%
%
%We then introduce an algorithm along with numerical results for different domain configurations. The algorithm is implemented using Fireshape, a Firedrake library that computes shape gradients via automatic differentiation and incorporates optimization algorithms from ROLTrilinos \cite{firedrake, fireshape, trilinos}.

\section{Problem statement}
\label{sec:problem_statement}
Let $\Omega \subset \mathbb{R}^d$ denote a bounded Lipschitz domain (with $d=2$ or $d=3$) representing a
neo-Hookean hyperelastic material with boundary $\partial \Omega$.
Dirichlet boundary conditions are imposed on $\Gamma_D \subset \partial \Omega$ with positive measure and natural traction-free boundary conditions are
imposed on $\Gamma_T \coloneqq \partial \Omega \setminus \Gamma_D$.
The displacement
\begin{equation}
    \label{eq:V}
    u\in V(\Omega)\coloneqq\left\{u\in H^1(\Omega,\R^d)\; :\; u|_{\Gamma_D}=u_D \right\}
\end{equation}
is a stationary point of the potential energy
\begin{equation}
     \E(u)\coloneqq\int_\Omega \psi(u) - b\cdot u \dx
\end{equation}
where \(b\in L^2(\Omega;\R^d)\) is the external force and $\psi(u)$ the internal energy density. Denoting by
$\F\coloneqq I+\nabla u$ the deformation gradient, the internal energy density $\psi(u)$ takes the form \cite[Section~6.4, pp.~182--188]{hyperelasticity_equation}
\begin{equation}
\label{eq:hyperelasticityEnergy}
\psi(u) \coloneqq  \frac{\mu}{2} \left( \text{tr}(\F^\top \F) - d \right) - \mu \ln\det(\F) + \frac{\lambda}{2} \left(\ln\det(\F)\right)^2\,,
\end{equation}
where \(\mu>0\) and \(\lambda>0\) are the Lamé coefficients. The Lamé coefficients are commonly expressed in terms
of the Young's modulus $E$ and the Poisson's ratio $\nu$ as $\mu=E/(2(1+\nu))$ and $\lambda=E\nu/((1+\nu)(1-2\nu))$. 
%Subsequently, the energy \(\E\) might also be referred  to as \(\E_\Omega\) for clarification.

By setting to zero the first order derivative of the potential energy $\E$, we obtain the nonlinear variational problem
\begin{equation}\label{eq:stateconstraint}
    \int_\Omega P(u):\nabla v \dx =\int_\Omega b\cdot v \dx \,\quad \text{for all } v\in V_0(\Omega)\,,
\end{equation}
where $V_0(\Omega) \coloneqq \{ u\in H^1(\Omega,\R^d)\,: u|_{\Gamma_D}=0\}$, and where 
\begin{equation}
\label{eq:PK}
P(u)\coloneqq \mu(\F-\F^{-\top})+\lambda \ln(\det(\F))\F^{-\top}
\end{equation}
denotes the first Piola--Kirchhoff stress tensor (see \Cref{appendix:PL}). 
Note that \eqref{eq:stateconstraint} is also the weak formulation of the nonlinear boundary value problem
\begin{equation}
    \label{eq:pde}
    \begin{aligned}
        -\nabla\cdot P(u)&=b && \text{in }\Omega, \\
        u&=u_D && \text{on }\Gamma_D, \\
        P(u)\cdot n&=0 && \text{on }\Gamma_T.\\
    \end{aligned}
\end{equation}

\begin{figure}[htbp]
    \centering
    \begin{subfigure}{0.48\textwidth}
        \centering
        \includegraphics[width=\linewidth]{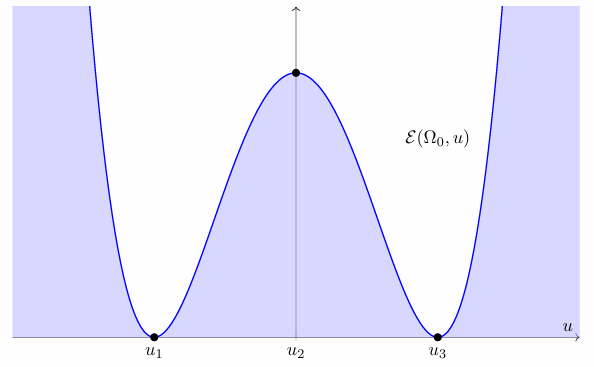}     
    \end{subfigure}
    \hfill
    \begin{subfigure}{0.48\textwidth}
        \centering
        \includegraphics[width=\linewidth]{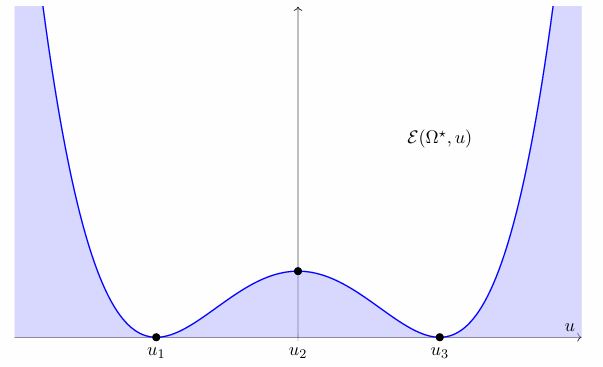}
    \end{subfigure}
    \caption{Energy jump before and after optimization.}
    \label{fig:energy-jump-before-and-after}
\end{figure}

Due to the nonlinearity of the Piola--Kirchhoff stress tensor $P$, the variational problem \cref{eq:stateconstraint}
generally admits multiple solutions. For a domain $\Omega$, let $u_1(\Omega),\dots,u_n(\Omega)$
denote the $n(\Omega) \in \N$ (distinct, isolated) solutions to \cref{eq:stateconstraint}.
To control the energy jump between solutions, we introduce the objective function
\begin{equation}
\label{eq:objective}
J(\Omega)\coloneqq\sum_{1 \le i < j \le n} C_{i,j}\left(\frac{\Delta\E_{i,j}(\Omega)}{\Delta\E_{i,j}(\Omega_0)} -r_{i,j} \right)^2\,,
\end{equation} where $\Omega_0$ denotes an initial domain, and where
\begin{equation}
\Delta\E_{i,j}(\Omega)\coloneqq \E\left(\Omega, u_i(\Omega)\right)-\E\left(\Omega, u_j(\Omega)\right)
\end{equation}
denotes the (signed) energy jump between the two solutions \(u_i(\Omega)\) and \(u_j(\Omega)\). 
% \footnote{Although this is the objective that we always plot and monitor, in the implementation we sometimes use the objective \begin{equation*}
% J(\Omega)\coloneqq\sum_{1 \le i < j \le n} C_{i,j}\left(\frac{1}{r_{i,j} }\frac{\Delta\E_{i,j}(\Omega)}{\Delta\E_{i,j}(\Omega_0)} -1\right)^2\,,
% \end{equation*} with appropriate weights $C_{i,j}$. This is because of specific implementation details. \ahs{The ratios $r_{i,j}$ need to be defined in the code as \texttt{self.r = firedrake.Constant(...)}, the reason is that for the continuation strategy (experiment 4), we need to dynamically update \texttt{r}, which is defined within the tape. The only reliable way I found for doing it is to use \texttt{firedrake.Constant}. Then the objective is computed as:
% \newline
% \smallskip
% {\ttfamily\small
% \noindent deltas = [\\
% \hspace*{4mm}assemble(\\
% \hspace*{8mm}(1 / self.r) * (self.energy(self.u[i]) - self.energy(self.u[j])) * dx,\\
% \hspace*{8mm}**self.assemble\_kwargs,\\
% \hspace*{4mm})\\
% \hspace*{4mm}for i, j in zip(self.objective\_params["i1"], self.objective\_params["i2"])\\
% ]\\
% \\
% obj = sum((delta / delta0 - 1) ** 2 for delta, delta0 in zip(deltas, self.deltas0))
% }
% \smallskip
% \newline
% where \texttt{delta0} is the initial energy jump (not divided by self.r).
% } }
This idea is illustrated in Figure \ref{fig:energy-jump-before-and-after}
where we have used the notation $\E(\Omega,u)$ to explicitly represent the dependency on the domain $\Omega$.

The coefficient $r_{i,j}\in\R$ encodes the desired jump for pair $(i, j)$, and the weight $C_{i,j}\geq 0$
balances the importance of different solution pairs (setting $C_{i,j}=0$ means that the jump between the solutions
 \(u_i(\Omega)\) and \(u_j(\Omega)\) is not controlled).
%\paragraph{Objective function}
%The PDE constraints introduced in the preceding section are non-linear, this usually leads to the existence of multiple solutions. For a domain \(\Omega\in \U(\Omega_0)\), let \(X(\Omega)=(u_1(\Omega),\cdots,u_d(\Omega))\) be a tuple of solutions to the PDE constraints (with \(d\in\mathbb{N}^*\) chosen), for now we allow repeating solutions so that \(d\) can be chosen arbitrarily (i.e if \(d\) is greater than the number of solutions we simply allow repetition). We define the energy jump
%\[\Delta_{i,j}(\Omega)=|\E(u_i(\Omega))-\E(u_j(\Omega))|\]
%for any two solutions \(u_i\) and \(u_j\) part of \(X(\Omega)\).
%
%Take  \(I_1, I_2\subset[[1,d]]\) two sets of indices. Our goal is to reduce the energy jump between the solutions, given by \(I_1\) and \(I_2\), to a ratio \(r\in(0,1)\) of the energy jump for an initial domain \(\Omega_0\). We therefore define the objective function
%\[J(\Omega)=\sum_{i\in I_1, j\in I_2} C_{i,j}\left(\Delta_{i,j}(\Omega) -r\Delta_{i,j}(\Omega_0) \right)^2\]
%where \( C_{i,j}\) are normalizing coeficients. 
Finally, to prevent distinct solutions from collapsing into each other, we consider the penalty term
\[B(\Omega)\coloneqq p\sum_{1 \le i < j \le n} \frac{C_{i,j}}{||u_i(\Omega)-u_j(\Omega)||_{L^2(\Omega)}^2}\]
where $p>0$ is a penalization parameter.

The associated shape optimization problem can be formulated as follows:
%find a domain $\Omega^*$ (among a set of admissible domains $\U$)
%that minimizes the objective
%\cref{eq:objective} constrained to the state equation \cref{eq:stateconstraint}, that is,
%\begin{equation}
%\label{eq:optimizationproblem}
%J(\Omega^*) =  \min_{\Omega\in\U} J(\Omega) +P(\Omega)\quad \text{subject to \cref{eq:stateconstraint}}\,.
%\end{equation}
\begin{equation}
\label{eq:optimizationproblem}
\begin{aligned}
\min_{\Omega \in \U(\Omega_0)} \quad & J(\Omega) + B(\Omega) 
\end{aligned}
\end{equation}
where the dependency of the functionals on $\Omega$ is through the integrals and implicitly through $u_i(\Omega)$, which are precisely the solutions of \eqref{eq:stateconstraint}.
Here $\U(\Omega_0)$ denotes a set of admissible domains constructed from a reference domain $\Omega_0$; more details are given in \Cref{sec:shapetheory} below.

Unlike most shape optimization problems, the functions \(u_i(\Omega)\)
appearing in the objective \cref{eq:objective} are all solutions to the same PDE constraint.
Since different domains may admit different numbers of solutions, we tacitly assume that
\cref{eq:optimizationproblem} is solved on the initial domain $\Omega_0\in\U$, where
all solutions are computed and properly labelled. This initial domain is then updated iteratively, and
at each domain update the state solutions are recomputed and properly labelled to ensure consistency.
This requires the existence of suitably regular maps $\Omega\mapsto u_i(\Omega)$
for each tracked solution \(u_i(\Omega)\). The existence of such maps is discussed in
\Cref{sec:shapetheory}.

\section{Iterative shape optimization}\label{sec:shapetheory}

For an initial domain \(\Omega_0\), the set of admissible controls is given by
\begin{equation}
    \label{eq:Uad}
    \U(\Omega_0)\coloneqq\left\{ \Omega=T(\Omega_0) \;:\; T\in \mathbf{F}(\Theta)  \right\}    
\end{equation}
where the set of admissible deformations is given by
\begin{equation}
\mathbf{F}(\Theta) \coloneqq \left\{I + \theta \;:\; \theta \in \Theta\,, 
\theta = 0\text{ on }\Gamma_D\,, (I + \theta)^{-1} \text{ exists, and } (I + \theta)^{-1} - I \in \Theta\right\} 
\quad
\end{equation}
for the Banach space \(\Theta = W^{1,\infty}(\mathbb{R}^d, \mathbb{R}^d)\). This makes \(\U(\Omega_0)\) a complete metric space \cite{delfour2011shapes}. Notice in particular that reachable shapes have the same topology as the reference domain. %In practice, \(\Omega_0\) serves as the base for constructing $\mathcal{U}$, the initial guess for the shape optimization, and the domain to which the PDE \eqref{eq:stateconstraint} is pulloften confounded with the initial guess for the shape optimization.

To solve the shape optimization problem \cref{eq:optimizationproblem}, we construct a sequence of diffeomorphisms
$\{T_k\}_{k\geq 1}$ of the form
\begin{equation}\label{eq:Tupdate}
T_k(x) = x + dT_k(x) \quad \text{for }x\in \Omega_0\,.
\end{equation}
Under suitable selection of the perturbation terms $\{dT_k\}_{k\geq 1}$, the
diffeomorphism sequence $\{T_k\}_{k\geq 1}$ induces a sequence of domains
$\{\Omega_k\coloneqq T_k(\Omega_0)\}_{k\geq 1}$ such that the
sequence $\{J(\Omega_k)\}_{k\geq 1}$ is decreasing. Our focus is on the
numerical realization of a shape optimization algorithm, \Cref{alg:optimization_loop}, to tackle \cref{eq:optimizationproblem}.
For this reason, we do not address the theoretical questions of existence and uniqueness
of solutions to \cref{eq:optimizationproblem} or determining conditions that guarantee
that $\{J(\Omega_k)\}_{k\geq 1}$ is a minimizing sequence.

The realization of optimization algorithms based on \cref{eq:Tupdate}
requires a guarantee that the solutions $u_i$ to the state constraint \cref{eq:stateconstraint}
can be consistently indexed and paired with the corresponding coefficients $r_{i,j}$ and weights
$C_{i,j}$ to evaluate \cref{eq:objective}. Under suitable assumptions, the IFT
provides this guarantee, as we now explain.

For a fixed vector field $A\in W^{1,\infty}(\Omega_0,\mathbb{R}^d)$ with $A=0$ on $\Gamma_D$, 
consider the family of transformations
$\{S_t(x) \coloneqq x +t\,A(x)\}_{t\geq0}$. It is well known that $S_t$ is a $W^{1,\infty}$-diffeomorphism if
the parameter $t$ is sufficiently small, say $t\leq\epsilon$ for an $\epsilon>0$ \cite{delfour2011shapes}.
In light of this, we introduce the function $G:[0,\epsilon]\times V(\Omega_0)\to V^{\star}(\Omega_0)$
defined as
\begin{equation}\label{eq:pulledback_form}
\langle G(t, u), v\rangle \coloneqq \int_{\Omega_0} \bigg( P_t(u)\colon (DS_t^{-\top}\nabla v)  - (b\circ S_t)\cdot v
\bigg)\det(DS_t)\dx\,,
\end{equation}
where $V^{\star}(\Omega_0)$ denotes the dual of $V_0(\Omega_0)$,
$\langle\,,\,\rangle$ denotes the $V^{\star}(\Omega_0)\times V_0(\Omega_0)$ duality pairing,
\begin{equation}
P_t(u)\coloneqq \mu(\F_t-\F_t^{-\top})+\lambda \ln(\det(\F_t))\F_t^{-\top}\,,
\end{equation}
and $\F_t\coloneqq I+DS_t^{-\top}\nabla u$, where $DS_t$ denotes the Jacobian of $S_t$.
Formula \cref{eq:pulledback_form} is obtained by pulling back
the weak formulation \cref{eq:stateconstraint} onto $\Omega_0$.
Therefore, the solutions to
\begin{equation}\label{eq_pullback_state_equation}
\langle G(t, u), v\rangle = 0 \quad \text{for all } v \in V_0(\Omega_0)\,,
\end{equation}
are the pullback of the solutions to \cref{eq:stateconstraint}, that is, if $u$ solves
\cref{eq_pullback_state_equation}, then the composition $u\circ S_t^{-1}$
solves \cref{eq:stateconstraint}, and vice-versa.
Moreover, if $\det(F_t)>0$, the function $G$ is continuous in $t$ and also continuously
differentiable in the variable $u$. Therefore,
if $(0,u)$ is a solution to \cref{eq_pullback_state_equation} and the Fr\'echet derivative $G_u(0,u)$
is invertible, the IFT \cite{IFT1,IFT2, IFT3,ciarlet2013}
guarantees the local existence of a differentiable function $t\to u(t)$ such that $G(t, u(t))=0$.
(This also implies that $u$ is an isolated solution.)
By assuming that $B_k\coloneqq T_{k+1}\circ T_{k}^{-1}-I\in W^{1,\infty}(\Omega_k)$ is sufficiently small
and applying the IFT to $S^k_t\coloneqq I + tB_k$, we can
verify that the solutions $u_i(\Omega_{k})$ and $u_i(\Omega_{k+1})$ correspond to the same
index $i$ by checking that the quantity
$\Vert u_i(\Omega_{k}) - u_i(\Omega_{k+1})\circ S^k_1\Vert_{L^2(\Omega_{k})}$
(or, equivalently but more conveniently,
$\Vert u_i(\Omega_{k})\circ (S^k_1)^{-1} - u_i(\Omega_{k+1})\Vert_{L^2(\Omega_{k+1})}$)
is also sufficiently small. Note that requiring $B_k$ to be
sufficiently small is not necessarily detrimental. Indeed, this guarantees that the pushed-forward
solution $u_i(\Omega_{k})\circ (S^k_1)^{-1}$ is a good initial guess for the computation of
$u_i(\Omega_{k+1})$ using an iterative nonlinear solver like Newton's method, which may thus
converge in just a few iterations.

To ensure that the sequence  $\{J(\Omega_k)\}_{k\geq 1}$ is decreasing,
we compute the perturbations $dT_{k}$
using the dogleg trust-region algorithm with BFGS-Hessian
updates implemented in the Rapid Optimization Library (ROL) \cite{rolsoftware}.
The underpinning steepest descent directions are based on
the shape gradients $\nabla J$ in $V_0(\Omega_0)$,
that is, the Riesz representative of the shape derivative of the constrained
function $J$ \cite{delfour2011shapes}, with respect to the
$V_0(\Omega_0)$ inner product \cite{paganini2018}. The whole optimization process is automated by
the shape optimization library Fireshape \cite{fireshape},
which builds on the finite element software Firedrake \cite{firedrake} and
leverages the automated adjoint derivation and (shape) differentiation of pyadjoint \cite{farrell2012c,pyadjoint} and of the
Unified Form Language (UFL) \cite{ufl, HaMiPaWe19}.

%\begin{algorithm}
%\caption{Optimization Loop \ap{This is not a trust-region algorithm}}
%\label{alg:optimization_loop}
%\begin{algorithmic}
%\REQUIRE Objective functional $J$, initial step size $\text{step}$, step size reduction factor $\alpha$, step size increase factor $\beta$, step tolerance $\text{step\_tol}$, maximum iterations $\text{max\_its}$, maximum decrease factor $\text{max\_decrease}$
%
%\STATE Initialize iteration counter $t \gets 0$
%\STATE Load solutions $(u_i(\Omega_0))$ from defcon analysis
%\STATE Compute initial objective value $J_{\Omega_0}$
%
%\WHILE{$t < \text{max\_its}$ \textbf{and} $\text{step} > \text{step\_tol}$}
%   \STATE Compute gradient $\nabla J(\Omega_t)$
%   \STATE Compute $\Tilde{\Omega} =( I - \text{step} \cdot \nabla J(\Omega_t))(\Omega_t)$
%   \STATE Compute solutions $(u_i(\Tilde{\Omega}))$ with Newton and initial guesses $(u_i(\Omega_t))$
%   \STATE Compute $J(\Tilde{\Omega})$
%
%   \IF{$J(\Tilde{\Omega}) < J(\Omega_t)$ \textbf{and} $J(\Tilde{\Omega}) \geq \text{max\_decrease} \cdot J(\Omega_t)$}
%       \STATE $t \gets t + 1$
%       \STATE $\Omega_t \gets \Tilde{\Omega}$
%       \STATE $(u_i(\Omega_t)) \gets $(u_i(\Tilde{\Omega}))$
%       \STATE $\text{step} \gets \text{step} \cdot \alpha$
%   \ELSE
%       \STATE $\text{step} \gets \text{step} \cdot \beta$
%   \ENDIF
%\ENDWHILE
%\end{algorithmic}
%\end{algorithm}

\begin{algorithm}
\caption{Shape optimization with trust region}
\label{alg:optimization_loop}
\begin{algorithmic}
\REQUIRE Objective functional $J$, initial domain $\Omega_0$, initial Hessian approximation $H_0$,
initial trust region radius $\Delta_0$, maximum radius $\Delta_{\max}$, step acceptance threshold $\eta_1 = 0.05$, radius shrinking threshold $\eta_1^- = 0.05$, radius growing threshold $\eta_2 = 0.9$, shrinking rates $\gamma_1^- = 0.0625$ (negative $\rho$), $\gamma_1^+ = 0.25$ (positive $\rho$), growth rate $\gamma_2=1.1$, sufficient decrease parameter $\mu = 10^{-4}$, gradient tolerance $\epsilon_g$, step tolerance $\epsilon_s$, maximum iterations $\text{max\_its}$
\STATE Initialize $k \gets 0$, $\Delta \gets \Delta_0$, $T_0\gets I$
\STATE Load initial solutions $(u_i(\Omega_0))$
\STATE Compute initial objective value $J(\Omega_0)$ and gradient $\nabla J(\Omega_0)$
\WHILE{$k < \text{max\_its}$ \textbf{and} $\|\nabla J(\Omega_k)\|_{V_0(\Omega_0)} > \epsilon_g$}
    \STATE Construct quadratic trust region model using L-BFGS Hessian approximation $H_k$:
    \begin{equation*}
        m_k(\theta) = J(\Omega_k) + \langle \nabla J(\Omega_k),\, \theta \rangle_{V_0(\Omega_0)}
        + \tfrac{1}{2}\langle H_k\,\theta,\, \theta \rangle_{V_0(\Omega_0)}
    \end{equation*}
    \STATE Solve trust region subproblem via Dogleg method:
    \begin{equation*}
        \theta_k = \arg\min_{\|\theta\| \leq \Delta} m_k(\theta)
    \end{equation*}
    \IF{$\|\theta_k\|_{V_0(\Omega_0)} \leq \epsilon_s$}
        \STATE \textbf{return}
    \ENDIF
    \STATE Compute candidate shape $\tilde{\Omega} = (T_k + \theta_k)(\Omega_0)$
    \STATE Compute solutions $(u_i(\tilde{\Omega}))$ with Newton and initial guesses $(u_i(\Omega_k)\circ T_k\circ (T_k + \theta_k)^{-1})$
    \STATE Compute $J(\tilde{\Omega})$ and actual-to-predicted reduction ratio:
    \begin{equation*}
        \rho_k = \frac{J(\Omega_k) - J(\tilde{\Omega})}{m_k(0) - m_k(\theta_k)}
    \end{equation*}
    \IF{$\rho_k \geq \eta_1$}
        \STATE $\Omega_{k+1} \gets \tilde{\Omega}$,  $T_{k+1} \gets T_k + \theta_k$, $(u_i(\Omega_{k+1})) \gets (u_i(\tilde{\Omega}))$
        \STATE Update gradient $\nabla J(\Omega_{k+1})$ and L-BFGS Hessian approximation $H_{k+1}$
        \STATE $k \gets k + 1$
    \ENDIF
    \STATE Update trust region radius $\Delta$ according to:
    \begin{equation*}
    \Delta \gets
    \begin{cases}
        \gamma_1^{-} \,\Delta & \text{if } \rho_k < 0 \\
        \gamma_1^{+} \,\Delta & \text{if } 0 \leq \rho_k < \eta_1^{-} \\
        \Delta              & \text{if } \eta_1^{-} \leq \rho_k < \eta_2 \\
        \min(\gamma_2\,\Delta,\, \Delta_{\max}) & \text{if } \rho_k \geq \eta_2 \text{ and } \Vert \theta_k\Vert = \Delta
    \end{cases}
    \end{equation*}
\ENDWHILE
\end{algorithmic}
\end{algorithm}

\section{Numerical Results}
\label{sec:numerical_results}
We demonstrate the effectiveness of the proposed methodology on a number of test cases in two and three dimensions.
In these numerical experiments, the Young's modulus is \(E=10^6\), the Poisson's ratio is \(\nu=0.3\),
and the internal body force is either \(b=(0,-10^3)^\top\) or \(b=(0,-10^3, 0)^\top\). The source code
for these experiments is available at \cite{codegithub}.

The initial computational domain $\Omega_0$
is an elastomeric square block with edges of length \(L=1\) and an embedded \(N\times N\)
square array of circular holes centered at \(\{(\frac{i}{N+1},\frac{j}{N+1})\,: i,j=1,\cdots,N\}\).
Horizontally (or vertically) adjacent holes have a center-to-center distance \(a = L/(N+1)\)
and a radius of \(R a/2\) with \(R\in(0,1)\).
The two vertical edges of $\Omega_0$ are flanked by a column of semicircles centered at
\(\{(x, i/(N+1))\,: i=1,\cdots,N, x=0,1\}\)
and of the same radius \(R a/2\).
Notice that if \(R=0\)  the circles are reduced to points while adjacent circle pairs touch each other
in a single point if \(R=1\). The parameter $R$ allows us to choose the complexity of the
multistable structure: we observe that the number of solutions increases as
\(R\) increases. Hereafter, we consider two cases:  \((N,R)=(2,0.9)\) and \((N,R)=(4,0.85)\).
%The fixed boundary $\Gamma_1$ is the lower edge of $\Omega^0$ and the prescribed displacement
%$u_0$ is applied on the upper edge of $\Omega_0$.
For three-dimensional simulations, $\Omega_0$ is further extruded in the third dimension
to a depth equal to \(h=1\).

The Dirichlet data is set to $u_D = 0$ on the bottom of the square, and $u_D=(0, -0.1)^\top$ or $u_D=(0, -0.1, 0)^\top$ on the top of the square. Natural boundary conditions are imposed on the remaining components of the boundary (side walls and walls of the holes).

To compute the set of initial solutions $\{u_i(\Omega_0)\}$, we use the deflated continuation algorithm \cite{farrell2026}.
The nonlinear equation \cref{eq:stateconstraint} is solved using PETSc's Newton algorithm, with the
linear systems solved with the sparse direct solver MUMPS \cite{petsc-web-page, petsc-user-ref, petsc-efficient, amestoy2001}.
Of the computed solutions, we select three specific solutions: two stable solutions $u_1(\Omega_0)$ and $u_3(\Omega_0)$ that are
$\mathbb{Z}_2$-symmetric with respect to each other and the saddle point $u_2(\Omega_0)$ in between.
We check the signs of the eigenvalues of the matrix corresponding to the discrete linearisation of \eqref{eq:stateconstraint}.
These are depicted in \cref{fig:Omega00} for a two-dimensional example. The approach below extends easily to other choices of solution branches.

%\begin{figure}[htb!]
%\includegraphics[scale=0.5]{Figures_paper/illustration-energy-landscape-before.png}
%\includegraphics[scale=0.5]{Figures_paper/illustration-energy-landscape-after.png}
%\caption{Sketch of the energy gap before (left) and after (right) optimization.}
%\end{figure}

We consider two possible families of objective functions:
% \pef{We already said what the objective function was, in \eqref{eq:objective}. It's different to the functionals below. Can we use that, or rewrite \eqref{eq:objective} to match what we used in experiments?}
\begin{align}
J_1(\Omega)=
\left(\frac{\Delta\E_{1,2}(\Omega)}
{\Delta\E_{1,2}(\Omega_0)}-r\right)^2
\quad\text{and}\quad
J_2(\Omega)=
\left(\frac{\Delta\E_{1,2}(\Omega)}
{\Delta\E_{1,2}(\Omega_0)}-r\right)^2
%\left((\frac{\E(u_1(\Omega))-\E(u_2(\Omega))}
%{\E(u_1(\Omega_0))-\E(u_2(\Omega_0))}-r\right)^2
+
\left(\frac{\Delta\E_{3,2}(\Omega)}
{\Delta\E_{3,2}(\Omega_0)}-r\right)^2\,.
%\left(\frac{\E(u_3(\Omega))-\E(u_2(\Omega))}
%{\E(u_3(\Omega_0))-\E(u_2(\Omega_0))}-r\right))^2
\end{align}
%The corresponding parameters are \ap{these are not consistent with \eqref{eq:objective},
%please fix them}
%\begin{align*}
%\text{either } &C_{21} = 1/\vert \Delta\E_{1,2}(\Omega_0)\vert\text{ and } C_{31} = C_{32} = 0\,,\quad\text{which define the objective function }J_1\\
%\text{or } &C_{21} = C_{32} = 1/\vert \Delta\E_{1,2}(\Omega_0)\vert\text{ and }C_{31} = 0\,,\quad\text{which define the objective function }J_2\,.
%\end{align*}
%this gives
%\begin{align*}
%J_1&=(\frac{\E_1-\E_2}{|\E_1-\E_2|}-r)^2\\
%J_1&=(\frac{\E_1-\E_2}{|\E_1-\E_2|}-r)^2+(\frac{\E_3-\E_2}{|\E_3-\E_2|}-r)^2
%\end{align*} 
%The other parameters are set to $r_{21} = r_{31} = r_{32} = r \vert \Delta\E_{1,2}(\Omega_0)\vert$,
%where $r=0.8, 0.6, 0.4, 0.2$ is a 
Here, $r$ denotes a target ratio.
Note that, in this particular case, the symmetry between $u_1(\Omega_0)$ and $u_3(\Omega_0)$ and
the central position of the saddle point $u_2(\Omega_0)$ (with respect to the energy functional $\E$) imply that
$\vert \Delta\E_{1,2}(\Omega_0)\vert = \vert \Delta\E_{2,3}(\Omega_0)\vert$.

\begin{figure}[!htb]
\centering
    \includegraphics[scale=0.15]{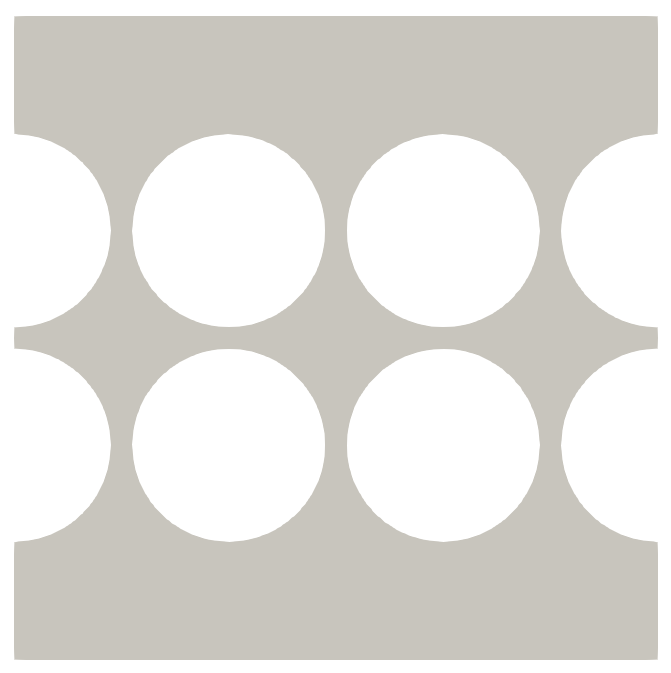}
    \includegraphics[scale=0.15]{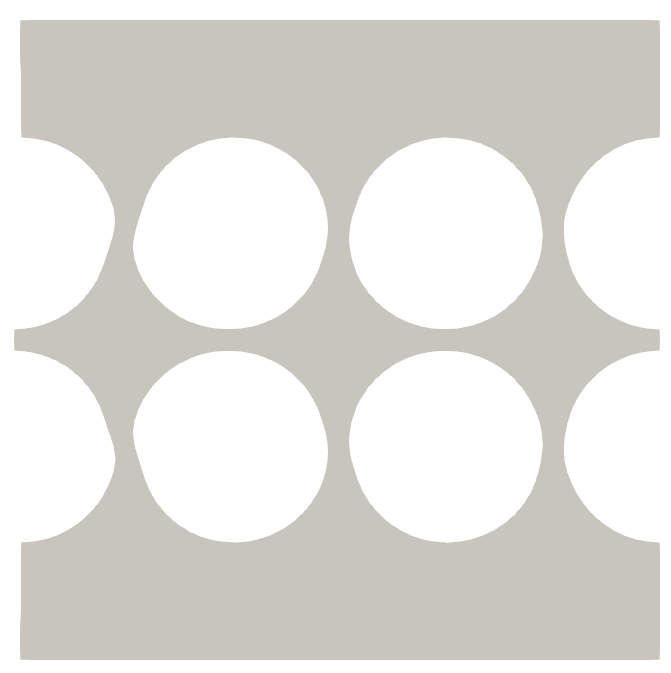} 
    \includegraphics[scale=0.15]{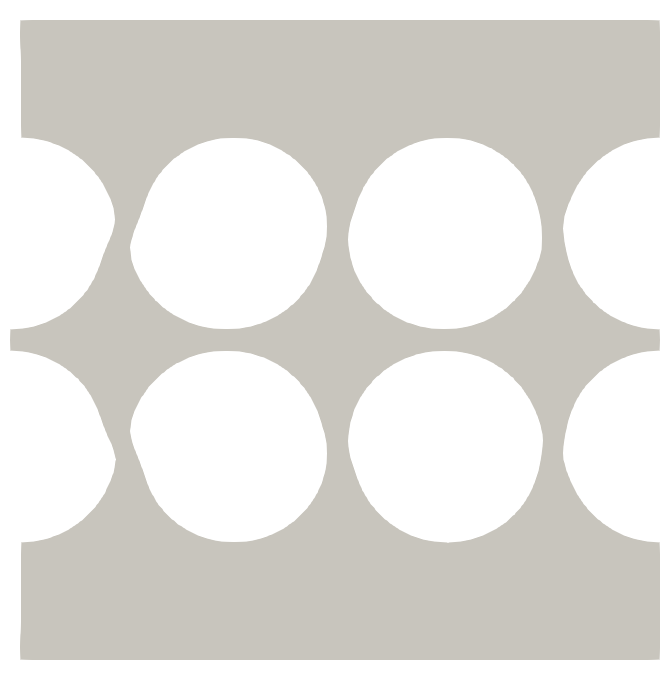}
    \includegraphics[scale=0.15]{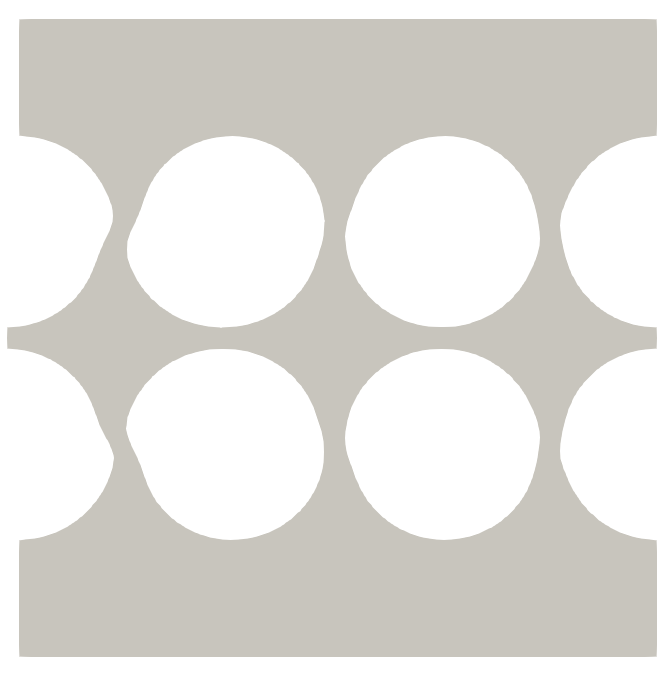}\\
    \includegraphics[scale=0.15]{Figures_paper/Omega0.png}
    \includegraphics[scale=0.15]{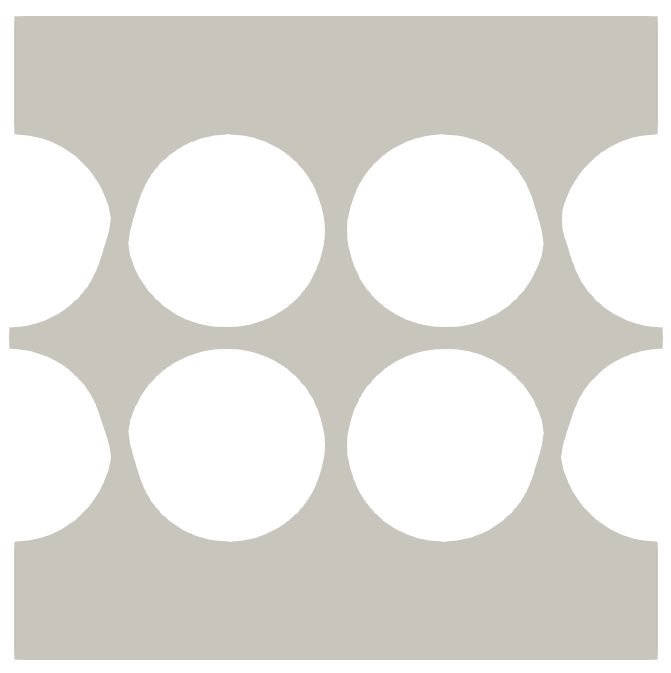}
    \includegraphics[scale=0.15]{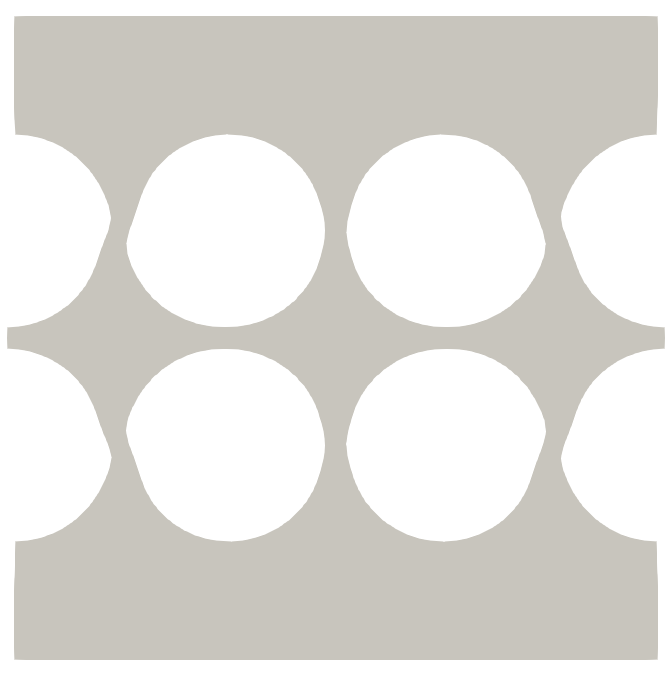}
    \includegraphics[scale=0.15]{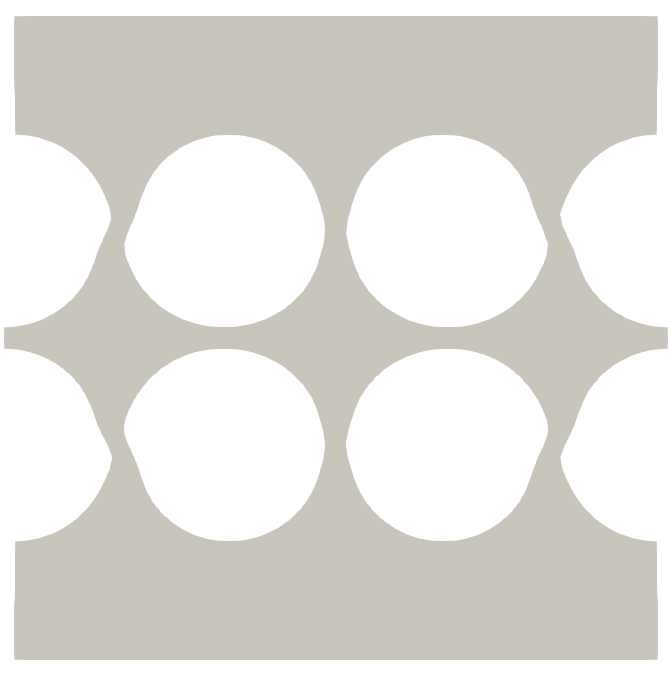}
\caption{Initial domain \(\Omega_0\) (leftmost), followed by the optimized shapes for the
target energy jump ratios \(r=0.6\), \(0.4\) and \(0.2\) (left to right),
with respect to \(J_1\) (first row) and \(J_2\) (second row).
%\ap{I suggest removing this plot. It takes up a lot of space and it is not very informative.}
%\pef{As I mention above, I think we should keep the $N=2$ and remove the $N=4$}
}
\label{fig:time_t}
\end{figure}

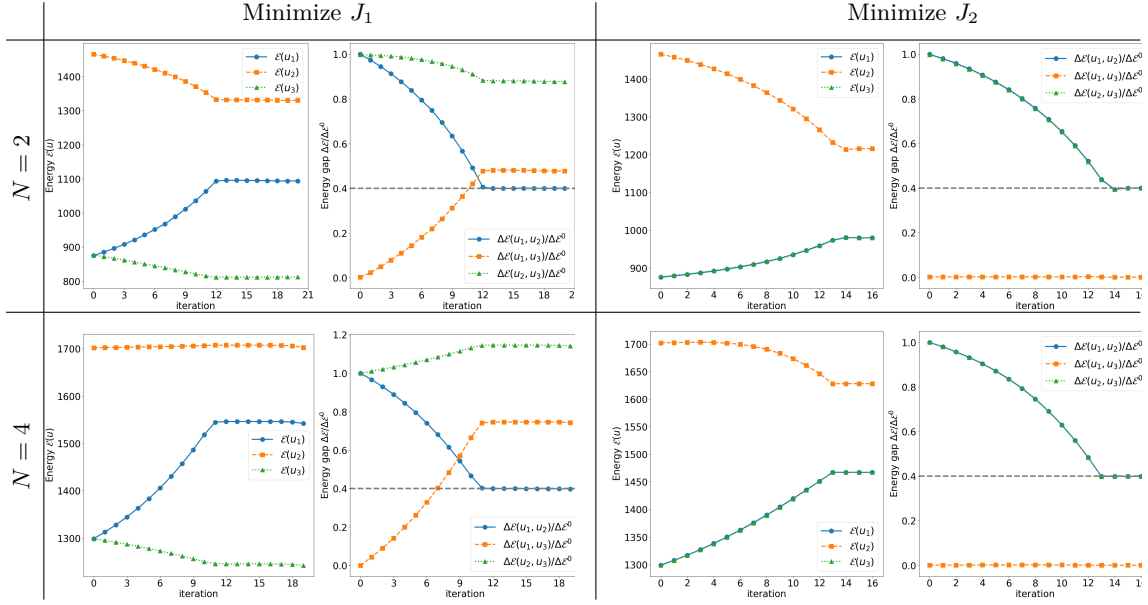
\begin{figure}[!htb]
\centering
\begin{tikzpicture}
    % Define styles for the row and column labels
    \tikzset{
        header/.style={font=\small},
        label_v/.style={header, rotate=90, anchor=center},
        label_h/.style={header, anchor=south}
    }

    % --- TOP BLOCK: N = 2 ---
    % Row Label
    \node[label_v] at (-0.3, 5.7) {$N=2$};
    
    % The Images (N=2)
    \node[anchor=south west, inner sep=0] (img1) at (0, 3.8) 
        {\includegraphics[scale=0.14, trim=0 0 740 0, clip]{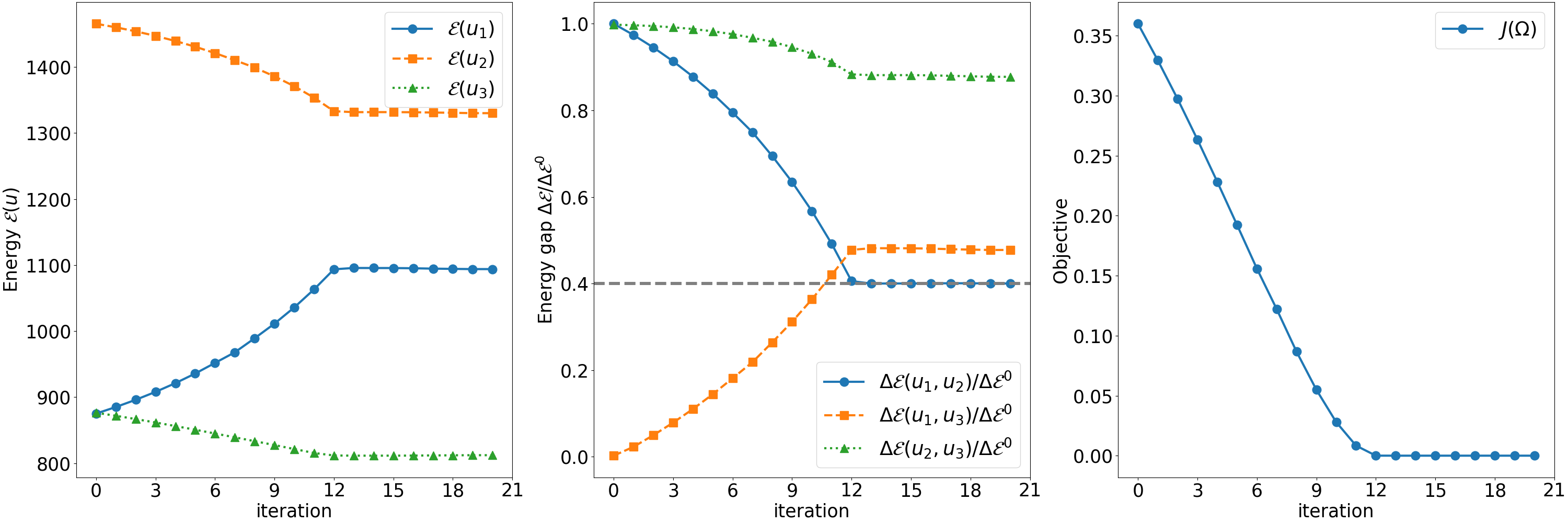}};
    \node[anchor=south west, inner sep=0] (img2) at (7.5, 3.8)
        {\includegraphics[scale=0.14, trim=0 0 740 0, clip]{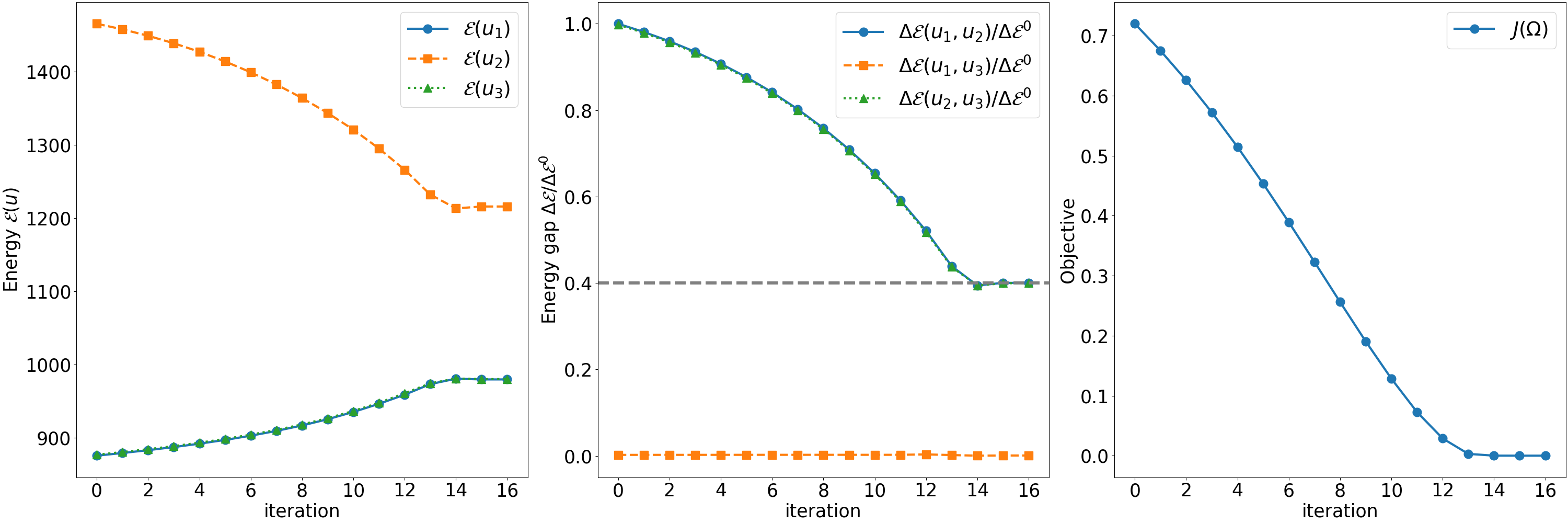}};

    % --- BOTTOM BLOCK: N = 4 ---
    % Row Label
    \node[label_v] at (-0.3, 1.9) {$N=4$};
    
    % The Images (N=4)
    \node[anchor=south west, inner sep=0] (img3) at (0,0) 
        {\includegraphics[scale=0.14, trim=0 0 740 0, clip]{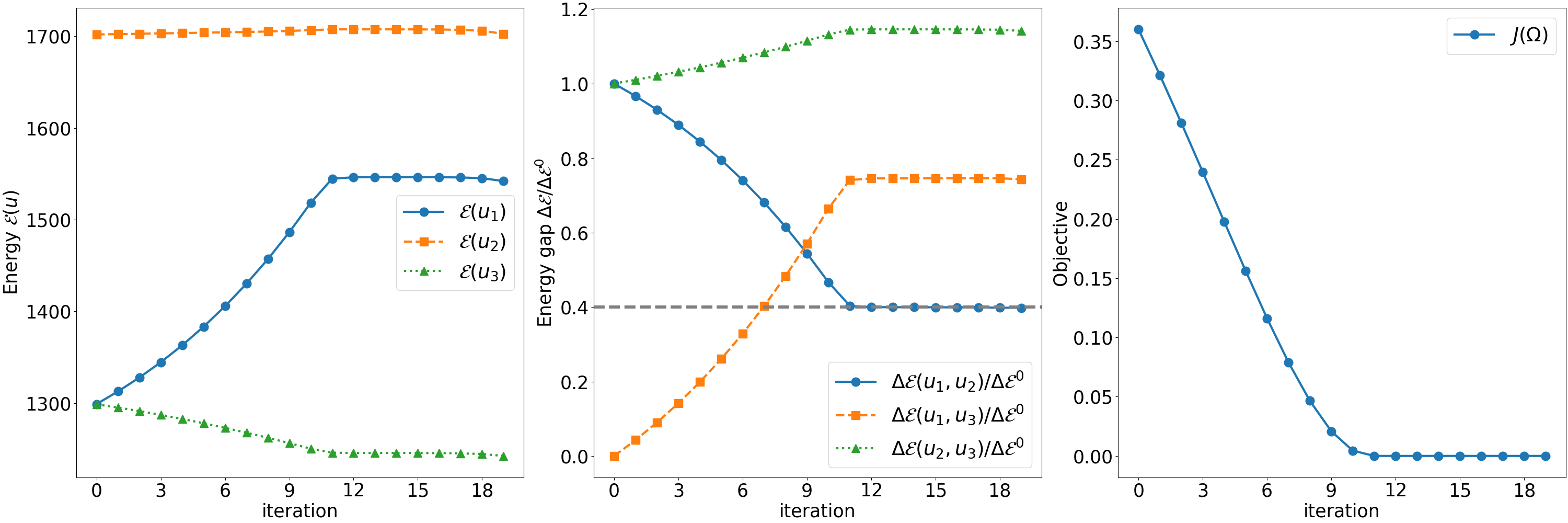}};
    \node[anchor=south west, inner sep=0] (img4) at (7.5,0)
        {\includegraphics[scale=0.14, trim=0 0 740 0, clip]{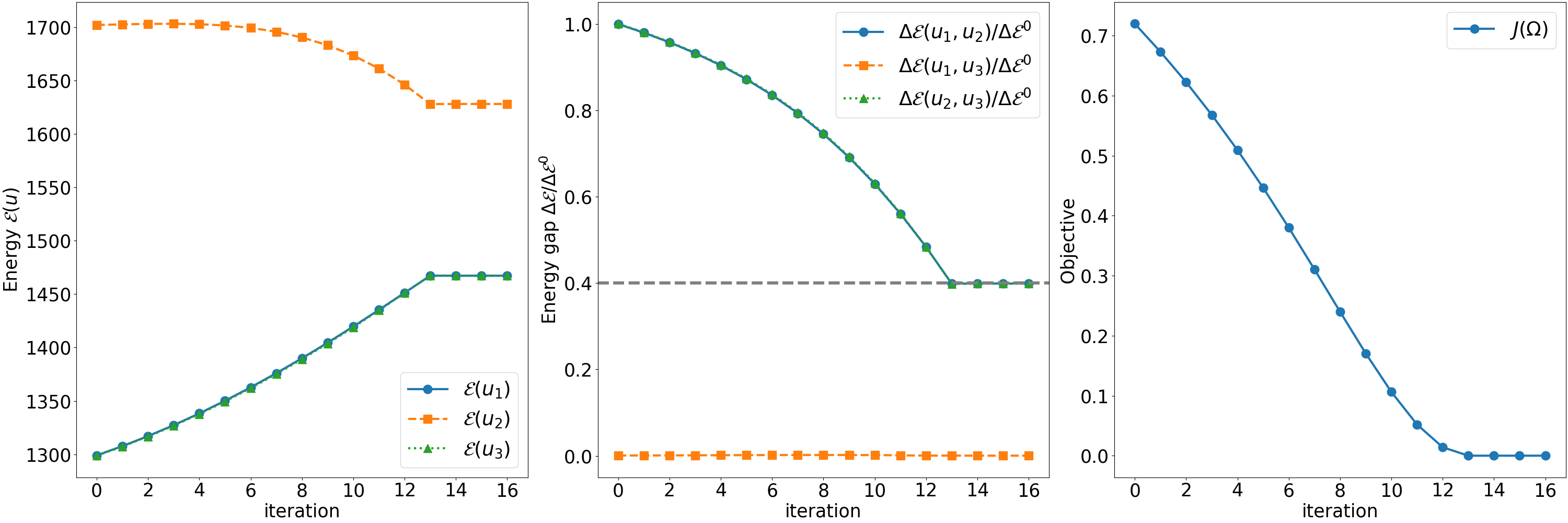}};

    % --- COLUMN HEADERS ---
    \node[label_h] at ( 3.5, 7.5) {Minimize $J_1$};
    \node[label_h] at (11.5, 7.5) {Minimize $J_2$};

    % --- GRID LINES (Optional) ---
    \draw (-0.5, 7.4) -- (14.5, 7.4);
    \draw (-0.5, 3.8) -- (14.5, 3.8);
    \draw (0, 0) -- (0,7.9);
     \draw (7.3, 0) -- (7.3,7.9);

\end{tikzpicture}

\caption{Convergence history for $r=0.4$ and $N=2,4$.
In the right column, the lines for $u_1$ and $u_3$ lie on top of one another.
}
\label{fig:2holes_NE2}
\end{figure}

\Cref{fig:time_t} illustrates, for the case $N=2$, the initial domain together with the optimized shapes obtained
for the target energy jump ratios \(r=0.6\), \(0.4\) and \(0.2\), using either objective \(J_1\) or \(J_2\).
In both instances, we observe nontrivial domain modifications. When the geometry is optimized by minimizing \(J_2\), which accounts for both gaps \(\Delta\E_{1,2}\) and \(\Delta\E_{2,3}\) simultaneously, we observe that the symmetry of the initial domain and its bistable structure are preserved throughout the shape optimization process.
In the convergence plots reported in \cref{fig:2holes_NE2}, we observe that when \(J_1\) is minimized, the energy curves of \(u_1\) and \(u_2\) converge toward each other, reflecting the targeted reduction of \(\Delta\E_{1,2}\)
to $40\%$ of $\Delta\E_{1,2}(\Omega_0)$ (dashed line).
In this simulation, the energy curve of \(u_3\), which is not controlled, drifts mildly.
In contrast, when minimizing \(J_2\), the energy curves of \(u_1\) and \(u_3\) coincide
(consistent with the observed symmetry preservation) and approach the energy curve of \(u_2\), reflecting the simultaneous
targeted reduction of \(\Delta\E_{1,2}\) and \(\Delta\E_{2,3}\).
\Cref{fig:2holes_NE2} also displays convergence plots for $N=4$, for which we observe a similar behavior.

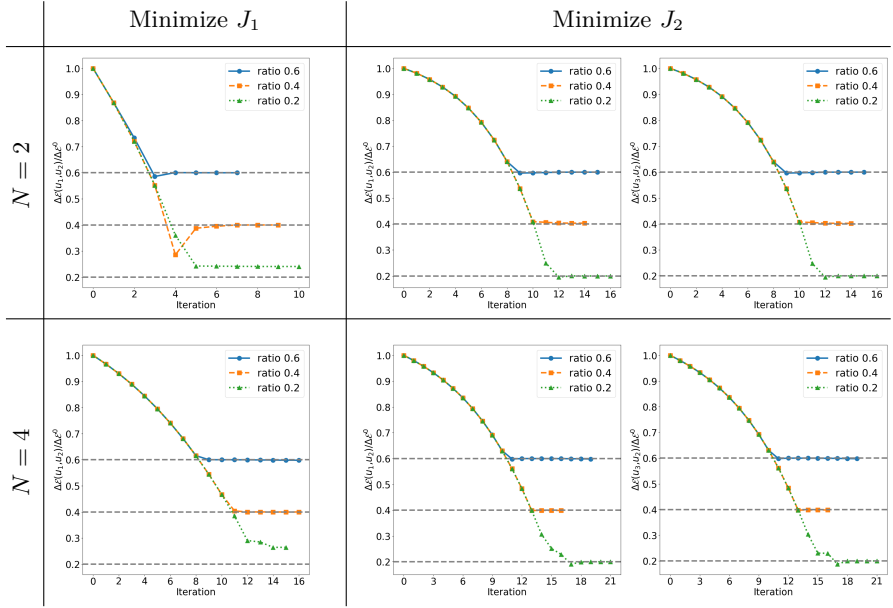
\begin{figure}[!htb]
\centering
\begin{tikzpicture}
    % Define styles for the row and column labels
    \tikzset{
        header/.style={font=\small},
        label_v/.style={header, rotate=90, anchor=center},
        label_h/.style={header, anchor=south}
    }

%	\begin{subfigure}[b]{0.8\textwidth}
%	\includegraphics[scale=0.15, trim=730 0 0 0, clip]{Figures_paper/ROL_NE3_2_holes_delta1.png}
%	\includegraphics[scale=0.15, trim=730 0 0 0, clip]{Figures_paper/ROL_NE3_2_holes_delta2.png}
%	\caption{Convergence history for \(N=2\) when minimizing \(J_1\) ({top}) and \(J_2\) ({bottom}).}
%	\end{subfigure}
%	
%	\begin{subfigure}[b]{0.8\textwidth}
%	\includegraphics[scale=0.15, trim=730 0 0 0, clip]{Figures_paper/ROL_NE3_4_holes_delta1.png}
%	\includegraphics[scale=0.15, trim=730 0 0 0, clip]{Figures_paper/ROL_NE3_4_holes_delta2.png}
%	\caption{Convergence history for \(N=4\) when minimizing \(J_1\) ({top}) and \(J_2\) ({bottom}). }
%	\end{subfigure}
	
  % --- TOP BLOCK: N = 2 ---
    % Row Label
    \node[label_v] at (-0.3, 5.7) {$N=2$};
    
    % The Images (N=2)
    \node[anchor=south west, inner sep=0] (img1) at (0.1, 3.8) 
        {\includegraphics[scale=0.14, trim=730 0 0 0, clip]{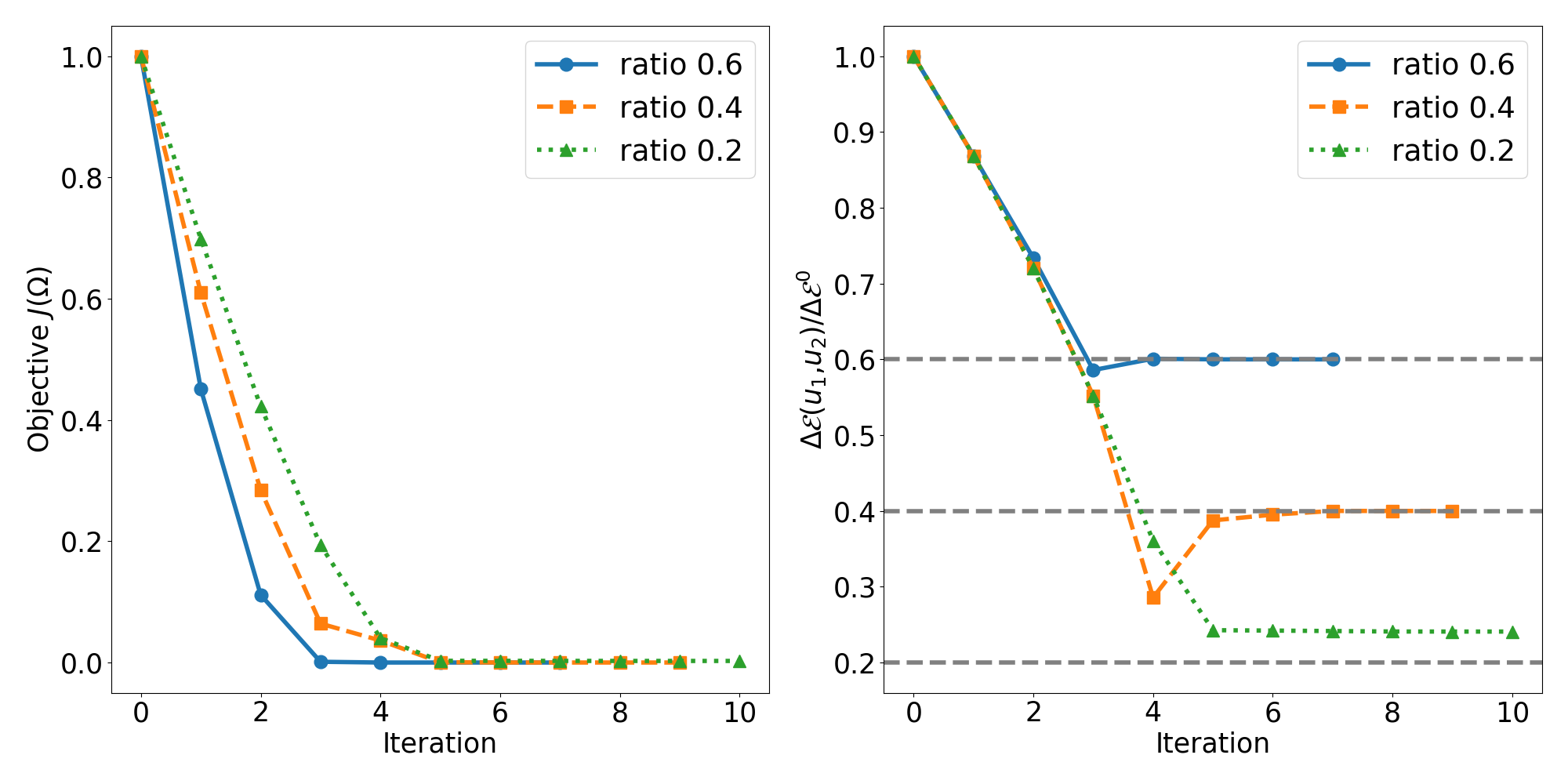}};
    \node[anchor=south west, inner sep=0] (img2) at (4.2, 3.8)
        {\includegraphics[scale=0.14, trim=730 0 0 0, clip]{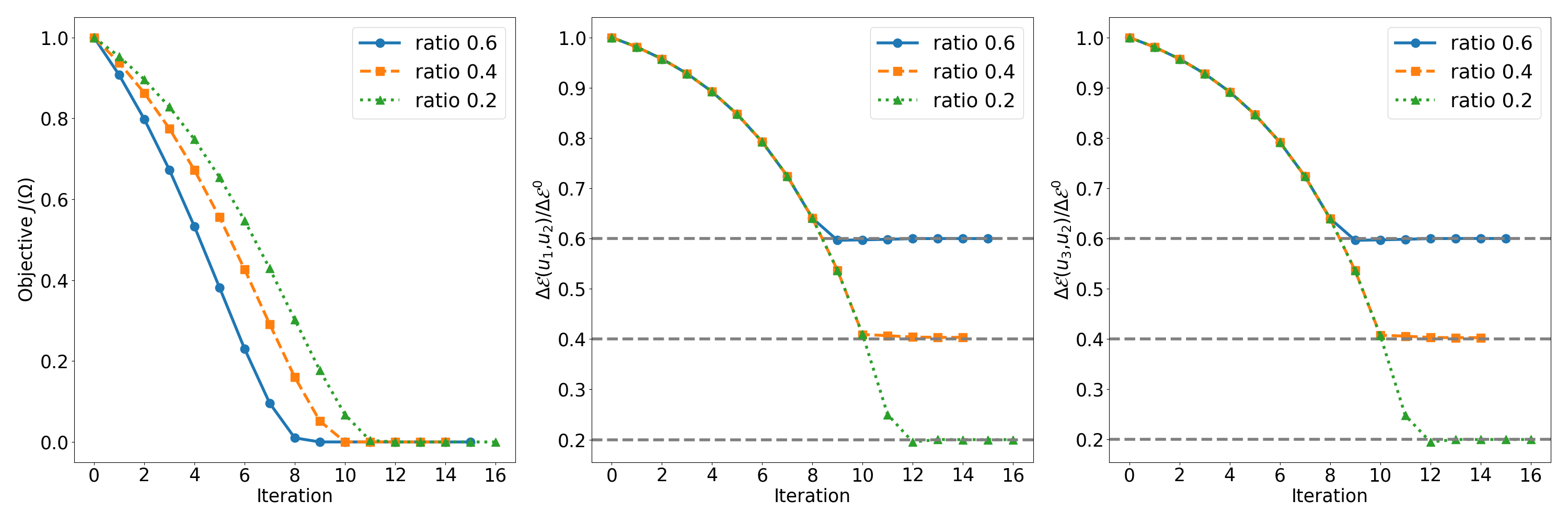}};

    % --- BOTTOM BLOCK: N = 4 ---
    % Row Label
    \node[label_v] at (-0.3, 1.9) {$N=4$};
    
    % The Images (N=4)
    \node[anchor=south west, inner sep=0] (img3) at (0.1,0) 
        {\includegraphics[scale=0.14, trim=730 0 0 0, clip]{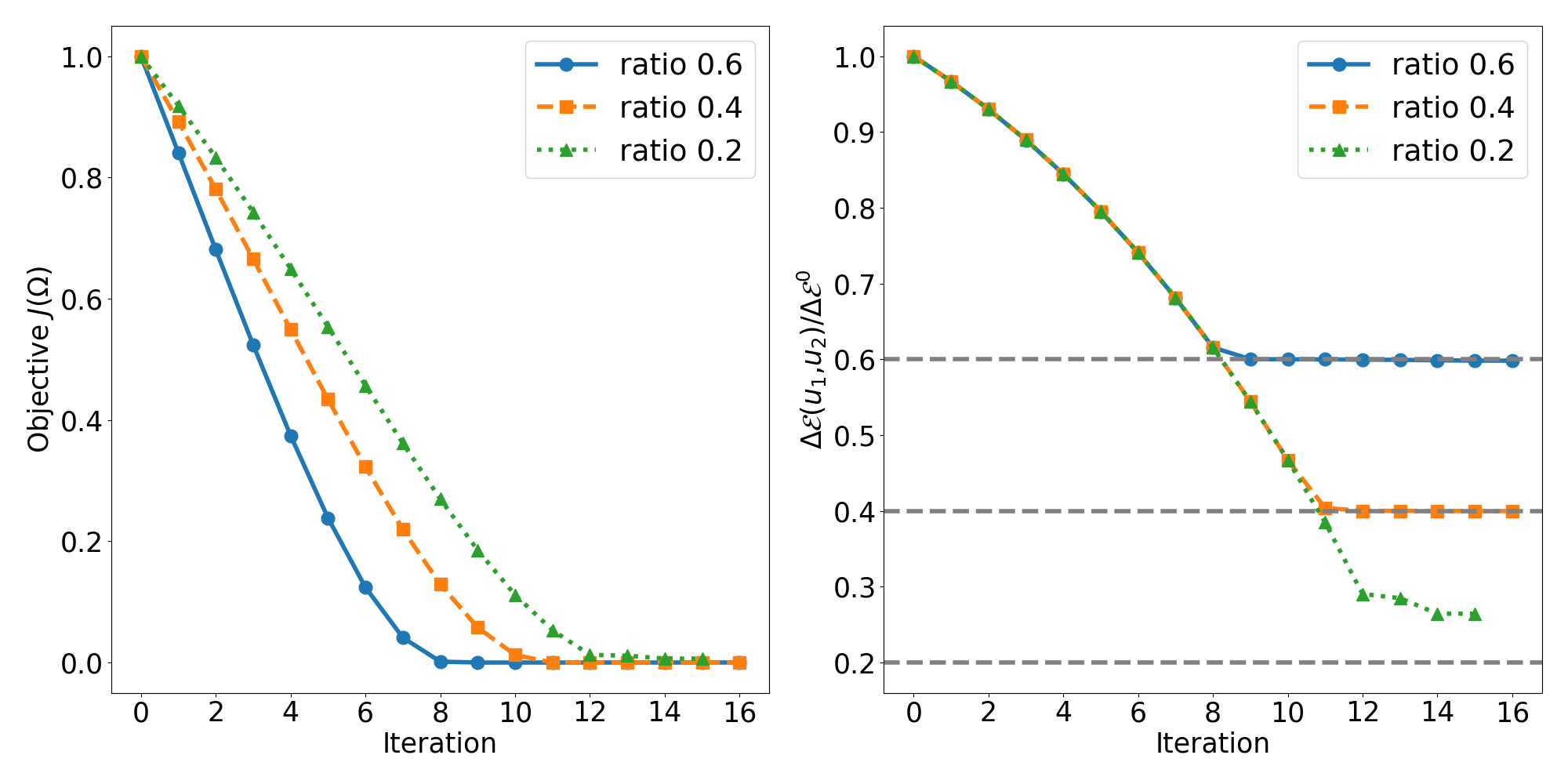}};
    \node[anchor=south west, inner sep=0] (img4) at (4.2,0)
        {\includegraphics[scale=0.14, trim=730 0 0 0, clip]{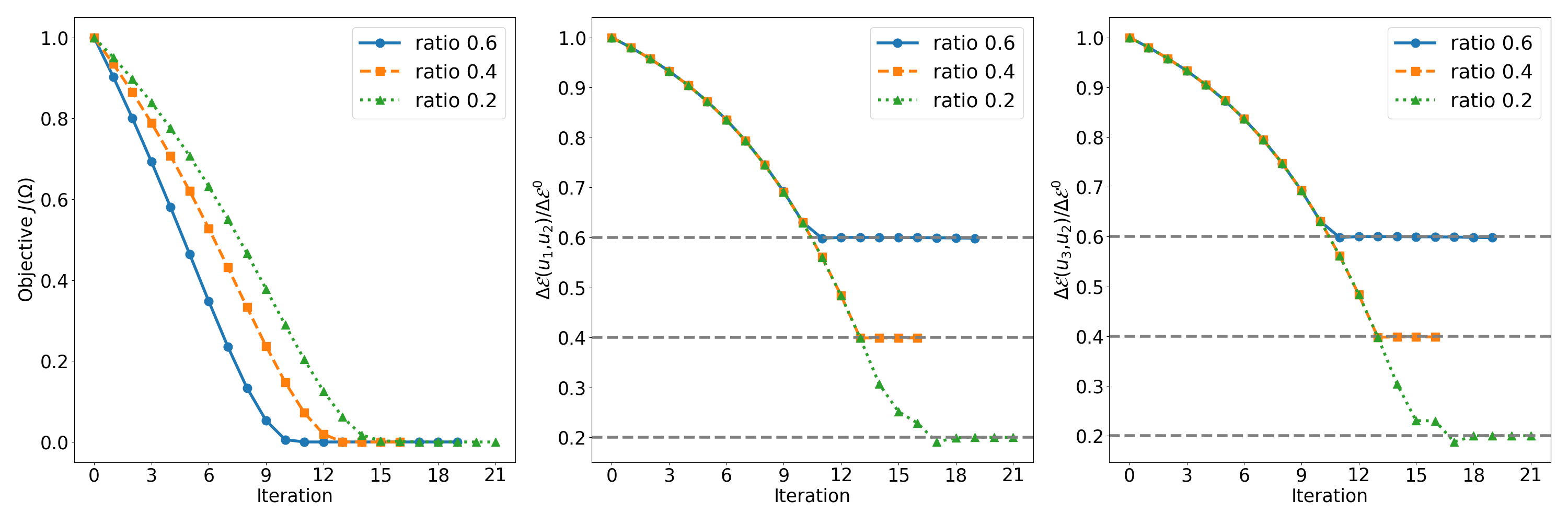}};

    % --- COLUMN HEADERS ---
    \node[label_h] at ( 2, 7.5) {Minimize $J_1$};
    \node[label_h] at (7.6, 7.5) {Minimize $J_2$};

    % --- GRID LINES (Optional) ---
    \draw (-0.5, 7.4) -- (11.2, 7.4);
    \draw (-0.5, 3.8) -- (11.2, 3.8);
    \draw (0, 0) -- (0,8);
     \draw (4, 0) -- (4,8);

\end{tikzpicture}
\caption{
Convergence history of energy gaps for \(r=0.2, 0.4, 0.6\) and $N=2, 4$. Dashed
lines indicate target ratios.  }
\label{fig:2holes_NE3}
\end{figure}

%\paragraph{Comparative convergence plots for different ratios}
We considered the same experiment, changing the target ratio ($r=0.2$, $0.4$, $0.6$). The corresponding
convergence histories are shown in \cref{fig:2holes_NE3}.
%shows convergence curves for the two objectives $J_{1}$ (top row) and $J_{2}$ (bottom row) for $N=2$ in panel (a) and $N=4$ in panel (b). 
%Each colored curve corresponds to a different target ratio $r$, representing the fraction of the initial energy gap ($\Delta\E_{12}(\Omega_0)$ or $\Delta\E_{23}(\Omega_0)$) to be retained at the optimum. The horizontal axis denotes iteration count, while the vertical axes display (left) the normalized objective and (middle and right) the individual energy gaps $\Delta\E_{12}$ and $\Delta\E_{23}$, respectively.
For smaller $r$ values, which demands more aggressive energy-gap reduction, the optimization problems take modestly more iterations.
In particular, the optimization for \(J_1\) and $r=0.2$ has more difficulty reaching the optimum and fails to reach the target ratio. Minimizing \(J_2\) for the same $r$ does achieve the target ratio.
We also consider the case \(r=1.6\), which corresponds to increasing the energy gap by \(60\%\). 
In this case, the desired ratio $r=1.6$ is achieved when minimizing
\(J_1\), but not when minimizing \(J_2\) (not shown).
We suspect that this different convergence behavior may be linked to the preservation (or loss) of symmetry,
but further research is needed
to confirm this observation. Comparing the results for $N=2$ and $N=4$, we observe that the optimization histories are qualitatively similar, indicating some robustness of the algorithm.

%\paragraph{Case $r>1$} 
%We also consider the case \(r=1.6\), which corresponds to increasing the energy gap by \(60\%\). 
%The optimized geometry for $N=4$ is shown in \cref{fig:OmegaT-r1.6}. In the geometry that minimizes $J_1$, we observe a slight stretch
%For the control of \(\Delta\E_{12}\) (top panel), the deformed geometry shows an adaptation where the domain slightly stretches 
%in alignment with the stable solution \(u_1\), effectively amplifying the energy barrier between solutions. The result is similar for $J_2$, in which case the result is symmetric and stretches
%Similarly, for the control of the two energy gaps \(\Delta\E_{12}\) and \(\Delta\E_{32}\) (bottom panel), the deformed geometry shows an adaptation where the domain slightly stretches
%in alignment with both stable solutions \(u_1\) and \(u_2\). 
%The corresponding convergence plots are portrayed in \cref{fig:r1.6-convergence}.
%In contrast to the case $r=0.2$, we observe that the desired ratio $r=1.6$ is achieved when minimizing
%\(J_1\), but not when minimizing \(J_2\).

\begin{figure}[!ht]
\centering
    \includegraphics[scale=0.2, trim=700 0 700 0, clip]{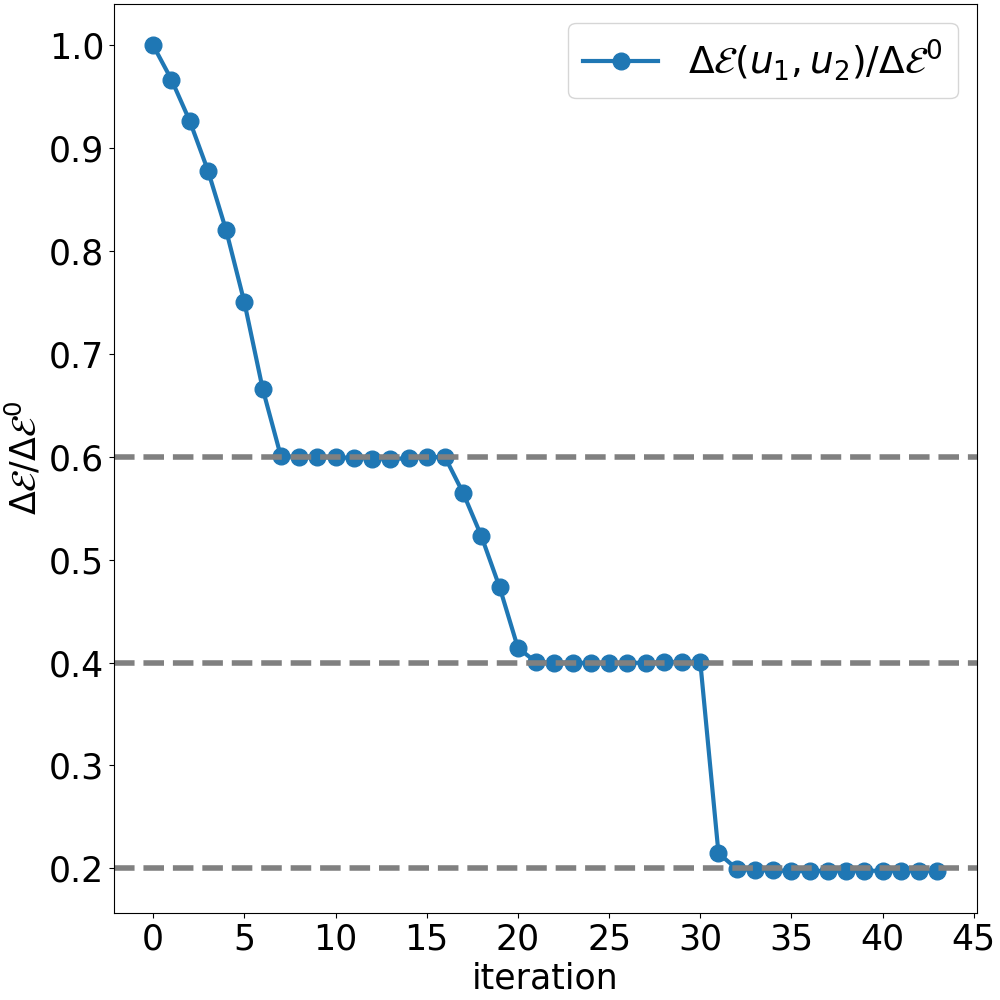}
\hspace{11cm}
    \includegraphics[scale=0.2, trim=700 0 700 0, clip]{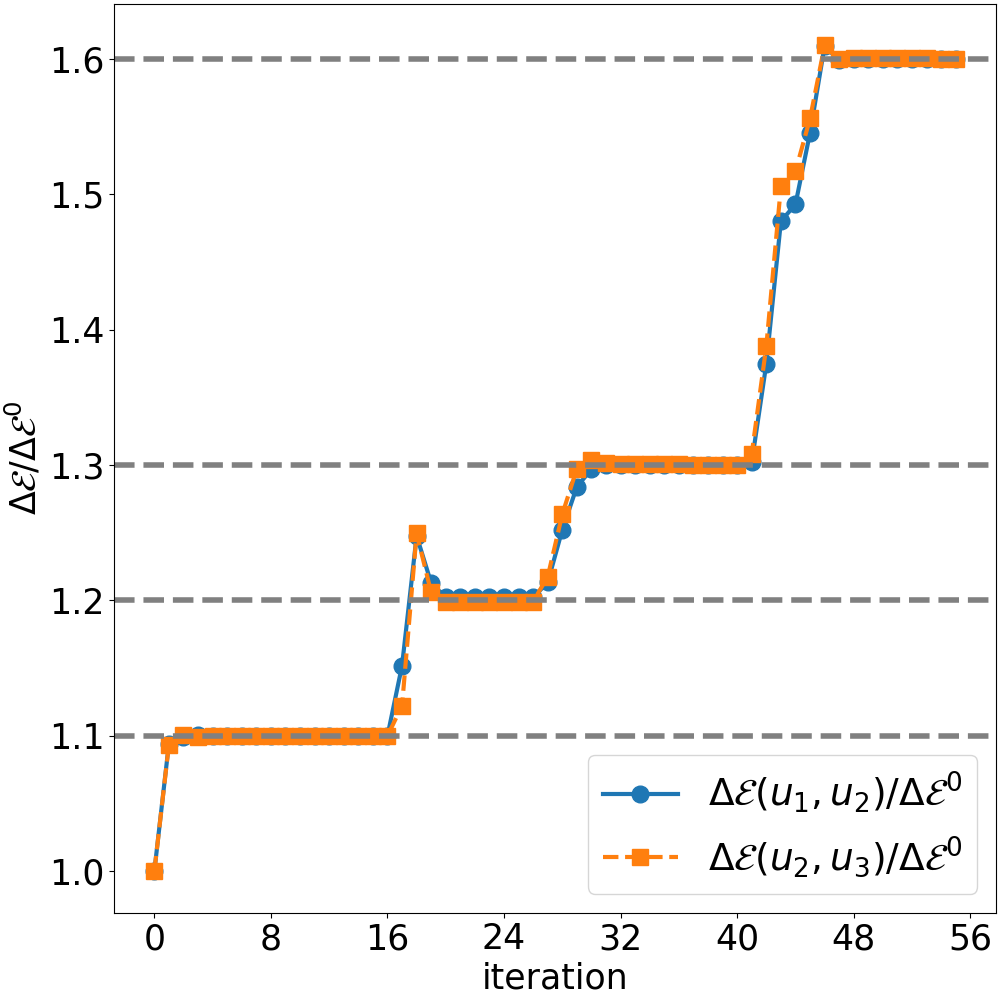}
\caption{\emph{Left:} sequential minimization of $J_1$ with target ratios $r=0.6, 0.4, 0.2$.
\emph{Right:} sequential minimization of $J_2$  with target ratios $r=1.1, 1.2, 1.3, 1.6$.
In both experiments, $N=4$. Dashed lines indicate intermediate target ratios.}
\label{fig:continuation}
\end{figure}

\begin{figure}
\centering
    \includegraphics[scale=0.2]{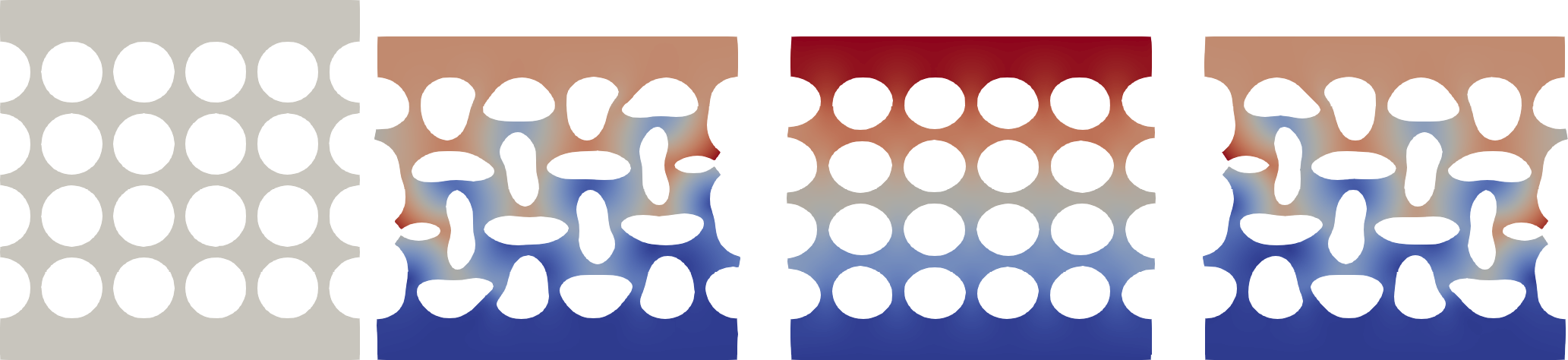}
    \\[0.5cm]
    \includegraphics[scale=0.2]{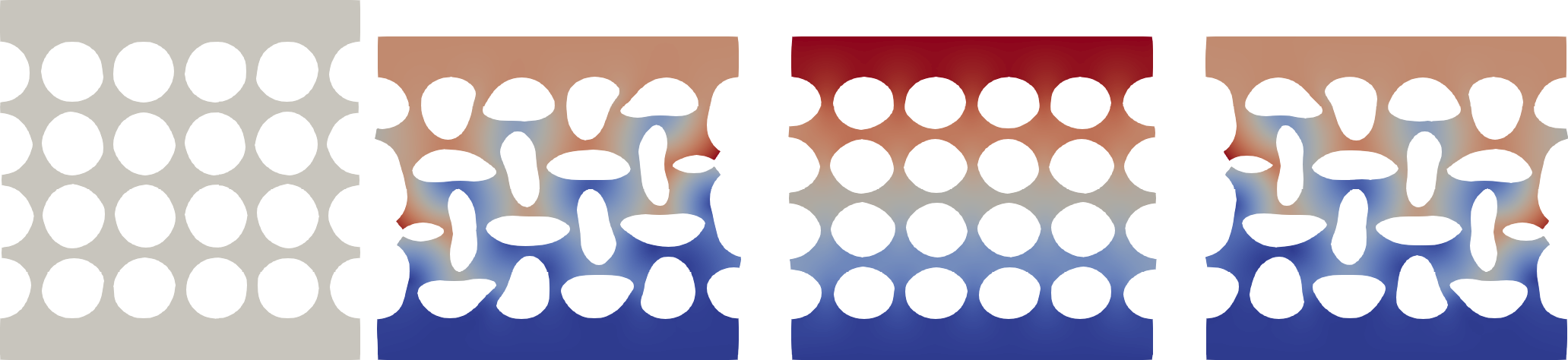}
    \\[0.5cm]
    \includegraphics[scale=0.2]{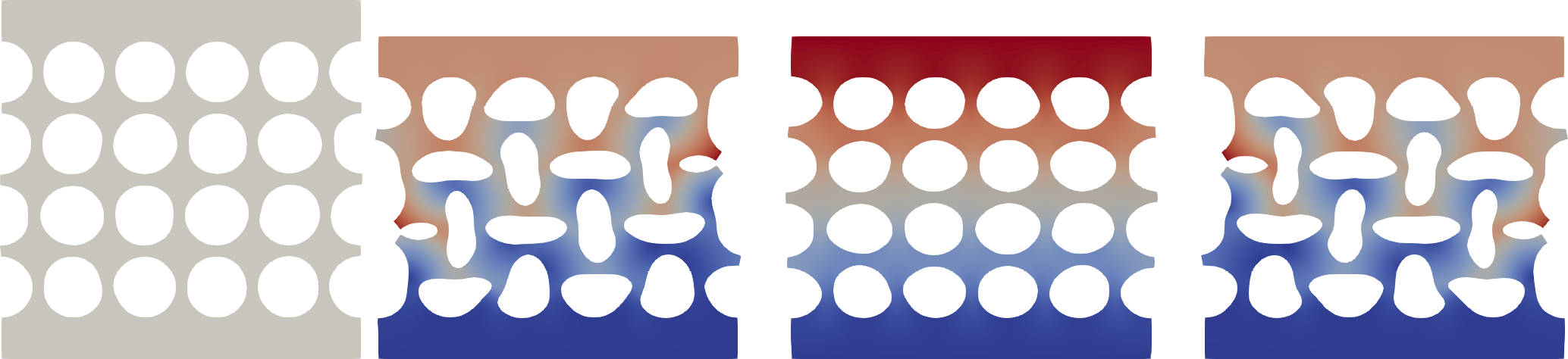}
    \caption{Domain $\Omega$ and its deformed configurations under displacements $u_1$, $u_2$ and $u_3$. From top to bottom: initial configuration, after optimization of $J_2$ with target ratio $r=1.6$, and after optimization of $J_1$ with target ratio $r=0.2$}
    \label{fig:continuation-geometries}
\end{figure}

To address the convergence failures observed in the previous experiments,
we consider a continuation strategy. Instead of directly targeting the ratio $r=0.2$ (or $r=1.6$),
we solve a sequence of shape optimization problems with decreasing (or increasing) target ratios $r$.
As evinced in \cref{fig:continuation} (for $N=4$),
this strategy successfully achieved the desired target ratios,
%we show the convergence history for \(J_1\)
%when the ratio is initially set to $r=0.6$ and, after convergence
%is reduced to $r=0.4$ and then to $r=0.2$. Similarly, in \cref{fig:continuation-r1.6},
%we consider the minimization of  \(J_2\) and increase the
%target ratio from $r=1.1$ to $r=1.6$ in five steps. In both case we set $N=4$.
%We observe a marked improvement, with the energy jumps achieving the desired ratio, 
albeit at the cost of an increased total number of iterations over the continuation steps. The initial and optimized geometries for this experiment are shown in \cref{fig:continuation-geometries}.
This demonstrates that while direct optimization for large energy gap increases can be challenging due to the competing geometric constraints, a carefully designed continuation strategy can overcome these difficulties and achieve the desired energy dissipation characteristics.

\begin{figure}[!ht]
\centering
\includegraphics[scale=0.25]{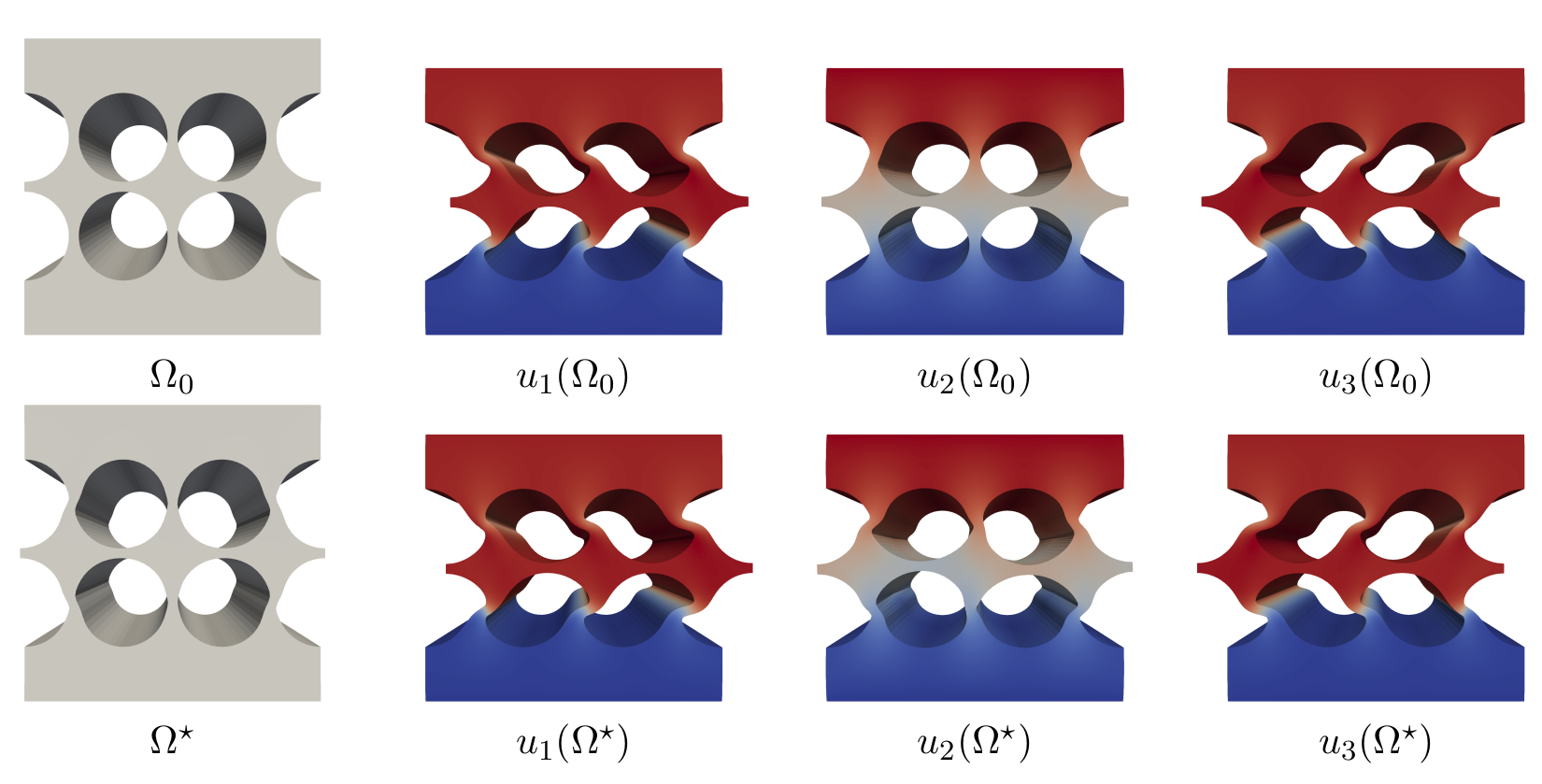}
%\\\includegraphics[scale=0.5]{Figures_paper/OptT.png}
\caption{Initial ({top}) and optimized ({bottom}) three-dimensional geometries for the target ratio \(r=0.2\) optimized by minimizing $J_2$. The leftmost geometry (in gray) denotes the optimized geometry at rest. From the left to the right, this geometry is displaced by applying  $u_1$, $u_2$, and $u_3$, respectively.}
\label{fig:3dN2}
\end{figure}

Finally, we perform a numerical experiment in three dimensions and minimize $J_2$ first
for $N=2$ and target ratio $r=0.2$. The initial and optimized domains (both at rest and
displaced) are shown in \cref{fig:3dN2}. %The corresponding convergence plots (not shown)
%are shown in \cref{fig:3D-convergence}.
%\cref{fig:3d} extends the analysis to a three-dimensional domain with $N=2$ and target ratio $r=0.2$. 
%The top row shows the initial geometry (gray) alongside its three solution branches (blue and red), while the bottom row presents the optimized shape and the associated solution branches. 
Consistent with the 2D case, the optimization successfully reduces the energy gap to approximately 20\% of its initial value (not shown), confirming the efficacy of the optimization strategy and objective formulation in higher dimensions. The optimized geometry redistributes material to moderate energy dissipation, with adaptations resembling the 2D pattern. In particular, the optimized domain is stretched opposite to the directions of the stable solution branches.
%In figure \cref{fig:3D-convergence}, we can see the convergence plots associated, we observe the same behaviour as in the 2 dimensional case.

\begin{figure}[!ht]
\centering
\includegraphics[scale=0.25]{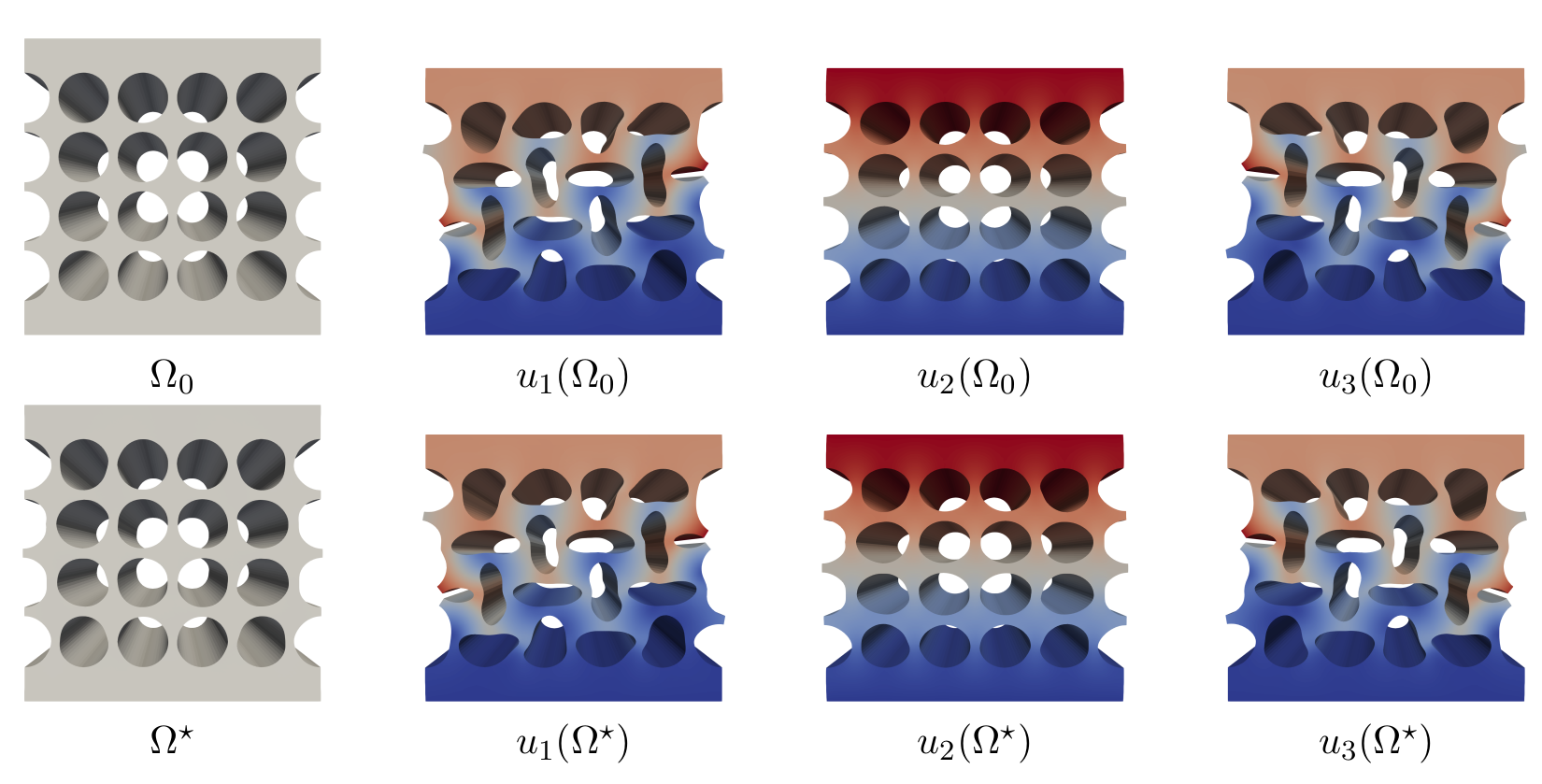}
%\\\includegraphics[scale=0.5]{Figures_paper/OptT.png}
\caption{Initial ({top}) and optimized ({bottom}) three-dimensional geometries for the target ratio \(r=0.2\) optimized by minimizing $J_2$ for $N=4$. The leftmost geometry (in gray) denotes the optimized geometry at rest. From the left to the right, this geometry is displaced by applying  $u_1$, $u_2$, and $u_3$, respectively.}
\label{fig:3dN4}
\end{figure}

Lastly, we reproduce the same experiment with $N=4$; the initial and optimized domains are shown in \cref{fig:3dN4}.

%\begin{figure}
%\centering
%\includegraphics[scale=0.15]{Figures_paper/NE2_N2_3D}
%%\caption{Comvergence plot for the 3 dimensional domain with \((N,r)=(2,0.9)\) for the optimization of \(J_2\) with \(r=0.2\). }
%\caption{
%Convergence history when minimizing $J_2$ with $N=2$ and
%target ratio $r=0.2$ for the three-dimensional case. 
%The plots portray the evolution of the displacement energies ({left}). of
%of the normalized energy gap ({center}), and of the objective function
%$J_2$ ({right}).}
%\label{fig:3D-convergence}
%\end{figure}
\section{Conclusion}
\label{sec:conclusion}
We have presented a shape optimization framework designed to control the energy jump between
distinct solutions of a neo-Hookean metamaterial, which has potential applications in impact energy dissipation.
The approach is general and can be applied to other shape optimization problems with 
multi-valued PDE constraints.

The methodology builds on rigorous theoretical foundations and its effectiveness
has been validated extensively through numerical experiments in two and three dimensions.

A central assumption of the proposed method is that the number of solution branches remains
constant throughout the optimization process. In scenarios where this does not hold,
a robust labelling procedure would be required to track branches of interest. Additionally, the current framework
cannot be directly applied to self-contact problems (such as those arising from large displacements),
in which case the PDE constraint is replaced by a quasi-variational inequality \cite{luo1996,hintermuller2011}.

\section*{Declaration} The research presented in this work was conducted by the authors. The use of AI tools was confined to assisting with text editing for spelling and grammar and preparing the code for reproducibility. The authors assume responsibility for all content.

\bibliographystyle{siam} 
\bibliography{biblio.bib}

\appendix
% Make \cref/\Cref say "Appendix A" instead of "Section A" after \appendix.
\crefalias{section}{appendix}
\crefalias{subsection}{appendix}
\crefalias{subsubsection}{appendix}
\section{Derivation of the first Piola--Kirchhoff stress tensor \eqref{eq:PK}}\label{appendix:PL}

Let $\Omega_0 \subset \R^3$ be a reference configuration.
A deformation $\varphi : \Omega_0 \to \R^3$ has deformation gradient
$\F = \nabla\varphi \in \R^{3\times 3}$, assumed to satisfy $J := \det\F > 0$.
We write $C = \F^\top\F$ for the right Cauchy--Green tensor and
$I_1 := \tr(C) = \tr(\F^\top\F)$ for its first invariant.

The Frobenius (double-dot) inner product of two matrices $A,B\in\R^{3\times 3}$
is defined by
\[
  A:B \;:=\; \tr(A^\top B)
                         \;=\; \sum_{i,j} A_{ij}\,B_{ij}.
\]
The first Piola--Kirchhoff (PK1) stress tensor $P\in\R^{3\times 3}$
is characterised by
\[
  \frac{\dd}{\dd\varepsilon}\bigg|_{\varepsilon=0} W(\F + \varepsilon\bH)
  \;=\; P : \bH
  \qquad\forall\,\bH\in\R^{3\times 3},
\]
i.e.\ $P$ is the matrix representative of the derivative of $W$ at $\F$.

We recall the neo-Hookean strain energy density \eqref{eq:hyperelasticityEnergy}
\begin{equation}\label{eq:W}
  W(\F) \;=\;
  \frac{\mu}{2}\bigl(\tr(\F^\top\F) - 3\bigr)
  \;-\; \mu\ln J
  \;+\; \frac{\lambda}{2}(\ln J)^2,
  \qquad J = \det\F,
\end{equation}
where $\mu > 0$ is the shear modulus and $\lambda \geq 0$ is the first
Lam\'{e} constant.
We split $W = W_1 + W_2 + W_3$ according to the three terms and compute
each derivative in turn.

\medskip
\noindent\textbf{Step 1---Linearisation of the trace term.}

Expanding $(\F+\varepsilon\bH)^\top(\F+\varepsilon\bH)$ gives
\[
  (\F+\varepsilon\bH)^\top(\F+\varepsilon\bH)
  \;=\; \F^\top\F
       + \varepsilon\bigl(\F^\top\bH + \bH^\top\F\bigr)
       + \varepsilon^2\,\bH^\top\bH.
\]
Taking the trace, the $O(\varepsilon^2)$ term vanishes upon differentiation at
$\varepsilon = 0$, yielding
\[
  \frac{\dd}{\dd\varepsilon}\bigg|_{\varepsilon=0}
  \tr\!\bigl((\F+\varepsilon\bH)^\top(\F+\varepsilon\bH)\bigr)
  \;=\;
  \tr(\F^\top\bH) + \tr(\bH^\top\F).
\]
Since $\tr(A) = \tr(A^\top)$ for any square matrix,
$\tr(\bH^\top\F) = \tr(\F^\top\bH)$, so
\[
  \frac{\dd}{\dd\varepsilon}\bigg|_{\varepsilon=0}
  W_1(\F+\varepsilon\bH)
  \;=\;
  \frac{\mu}{2}\cdot 2\,\tr(\F^\top\bH)
  \;=\;
  \mu\,\F:\bH.
\]
The contribution of $W_1$ to $P$ is therefore $\mu\F$.

\medskip
\noindent\textbf{Step 2---Linearisation of $\det$ via Jacobi's formula.}

Factor out $\F$ from the perturbed determinant:
\[
  \det(\F+\varepsilon\bH)
  \;=\; \det\F\cdot\det\!\bigl(\I + \varepsilon\F^{-1}\bH\bigr).
\]
The first-order Taylor expansion of the determinant near the identity reads
$\det(\I + \varepsilon A) = 1 + \varepsilon\,\tr(A) + O(\varepsilon^2)$,
so
\[
  \frac{\dd}{\dd\varepsilon}\bigg|_{\varepsilon=0}
  \det(\F+\varepsilon\bH)
  \;=\; J\,\tr(\F^{-1}\bH).
\]
Since $\tr(\F^{-1}\bH) = \tr\bigl((\F^{-\top})^\top\bH\bigr) = \F^{-\top}:\bH$,
this is Jacobi's formula:
\begin{equation}\label{eq:jacobi}
  \frac{\dd}{\dd\varepsilon}\bigg|_{\varepsilon=0}
  \det(\F+\varepsilon\bH)
  \;=\; J\,\F^{-\top}:\bH.
\end{equation}

\medskip
\noindent\textbf{Step 3---Linearisation of $\ln J$ and differentiation of $W_2$ and $W_3$.}

Applying the chain rule to $\ln$ and using \eqref{eq:jacobi},
\[
  \frac{\dd}{\dd\varepsilon}\bigg|_{\varepsilon=0}
  \ln\det(\F+\varepsilon\bH)
  \;=\;
  \frac{1}{J}
  \cdot J\,\F^{-\top}:\bH
  \;=\;
  \F^{-\top}:\bH.
\]

For $W_2 = -\mu\ln J$, this directly gives
\[
  \frac{\dd}{\dd\varepsilon}\bigg|_{\varepsilon=0}
  W_2(\F+\varepsilon\bH)
  \;=\;
  -\mu\,\F^{-\top}:\bH,
\]
contributing $-\mu\F^{-\top}$ to $P$.

For $W_3 = \tfrac{\lambda}{2}(\ln J)^2$, setting
$g(\varepsilon) := \ln\det(\F+\varepsilon\bH)$ and applying the chain rule to
$\tfrac{\lambda}{2}[g(\varepsilon)]^2$ at $\varepsilon = 0$,
\[
  \frac{\dd}{\dd\varepsilon}\bigg|_{\varepsilon=0}
  W_3(\F+\varepsilon\bH)
  \;=\;
  \lambda\,\ln J\cdot g'(0)
  \;=\;
  \lambda\ln J\;\F^{-\top}:\bH,
\]
contributing $\lambda\ln J\;\F^{-\top}$ to $P$.

\medskip
\noindent\textbf{Conclusion.}

Summing the three contributions, for every $\bH\in\R^{3\times 3}$,
\[
  \frac{\dd}{\dd\varepsilon}\bigg|_{\varepsilon=0} W(\F+\varepsilon\bH)
  \;=\;
  \bigl(\mu\F \;-\; \mu\F^{-\top} \;+\; \lambda\ln J\;\F^{-\top}\bigr):\bH.
\]
Since $P$ is defined by $P:\bH$ equalling this quantity for all $\bH$,
we conclude
\[
    P(\F)
    \;=\;
    \mu\bigl(\F - \F^{-\top}\bigr)
    \;+\;
    \lambda\ln(\det\F)\;\F^{-\top}.
\]

\end{document}